\documentclass{article}%
\usepackage{amsfonts}
\usepackage{amsmath}%
\usepackage{amssymb}%
\usepackage{graphicx}
\newtheorem{theorem}{Theorem}

\newtheorem{lemma}[theorem]{Lemma}

\newtheorem{proposition}[theorem]{Proposition}
\newtheorem{remark}[theorem]{Remark}

\newenvironment{proof}[1][Proof]{\noindent\textbf{#1.} }{\ \rule{0.5em}{0.5em}}
\begin{document}

\title{Three dimensional vortices in Abelian Gauge Theories}
\author{Vieri Benci$^{\ast}$, Donato Fortunato$^{\ast\ast}$\\$^{\ast}$Dipartimento di Matematica Applicata \textquotedblleft U.
Dini\textquotedblright\\Universit\`{a} degli Studi di Pisa Universit\`{a} di Pisa\\via Bonanno 25/b, 56126 Pisa, Italy\\e-mail: benci@dma.unipi.it\\$^{\ast\ast}$Dipartimento di Matematica \\Universit\`{a} di Bari and INFN sezione di Bari\\Via Orabona 4, 70125 Bari, Italy\\e-mail: fortunat@dm.uniba.it}
\date{}
\maketitle

\begin{abstract}
In this paper we consider an Abelian Gauge Theory in $\mathbb{R}^{4}$ equipped
with the Minkowski metric. This theory leads to a system of equations, the
Klein-Gordon-Maxwell equations, which provide models for the interaction
between the electromagnetic field and matter. A three dimensional vortex is a
finite energy solution of these equations in which the magnetic field looks
like the field created by a finite solenoid. Under suitable assumptions, we
prove the existence of vortex-solutions.

\end{abstract}
\tableofcontents

\section{Introduction\label{ded}}

Abelian gauge theories (in $\mathbb{R}^{4}$ equipped with the Minkowski
metric), provide models for the interaction between the electromagnetic field
and matter. Actually an Abelian gauge theory leads to a system of equations,
the Klein-Gordon-Maxwell equations (KGM), which occurs in various physical
problems (elementary particles, superconductivity, ...); see e.g.
\cite{rub},\cite{yangL}.

The properties of the solutions of the KGM depend on the lower order term $W$
(see (\ref{marisa})). Actually the choice of this term determines the
peculiarities of the various models.

In this paper we show that a suitable choice of the term $W$ guarantees the
existence of finite energy vortices in three space dimensions. Roughly
speaking, a vortex is a finite energy solution in which the magnetic field
looks like the field created by a solenoid.

As far as we know, the existence of vortex-solutions for Abelian gauge
theories has been studied only in the case of two space dimensions (see the
pioneering papers \cite{ab}, \cite{nil} and the books \cite{fel}, \cite{raj},
\cite{rub}, \cite{yangL} with their references). Clearly the two dimensional
vortices in the $x_{1}$, $x_{2}$ plane can be extended to $\mathbb{R}^{3}\;$as
constant maps in the $x_{3}$-direction. Of course these solutions have
infinite energy. We point out that, in the 2-dimensional models, the functions
$W$ that have been considered are of the type
\[
W(s)=\left(  1-s^{2}\right)  ^{2}%
\]
namely double well shaped and positive functions.

Here the lower order term which we have considered is the following one:
\begin{equation}
W(s)=\frac{1}{2}s^{2}-\frac{s^{p}}{p},\text{ }2<p<6,\text{ }s\geq0\label{oc}%
\end{equation}

The KGM with the lower order term (\ref{oc}) has been studied in \cite{bf},
\cite{befo}, \cite{ca}, \cite{tea}, \cite{tea2}, \cite{dav}. In these papers
the existence of stationary solutions having spherical symmetry has been
proved. It is easy to see that the spherical symmetry implies that both the
magnetic field and the angular momentum vanish. On the contrary, a
vortex-solution breaks the spherical symmetry (cf. Remark \ref{v}), moreover
the magnetic field and the angular momentum do non vanish.

\ 

Since the KGM are invariant for the Lorentz group, a Lorentz transformation of
a vortex creates a travelling \textit{solitary wave}. By \textit{solitary
wave}, we mean a solution of a field equation whose energy travels as a
localized packet. In this respect solitary waves have a particle-like behavior
(see e.g.\cite{fel}, \cite{raj}). A detailed discussion of solitary waves for
KGM can be found in \cite{befogranas}. Solitary waves obtained by vortices
behave as particles having a magnetic field and a kind of spin. When
considering evolution problems relative to KGM, the request $W(s)\geq0$ seems
necessary to have \textquotedblright good\textquotedblright\ solutions (cf.
\cite{befogranas}). In \cite{befo}, the existence of solitary waves relative
to stationary solutions having spherical symmetry has been proved also for a
class of functions $W\geq0.$ The existence of three dimensional vortices for
positive $W$ is still an open problem.

The paper is organized as follows: In the first section we introduce the
KGM-equations, we give a precise definition of three dimensional vortex
solution and we state the main existence theorem. In the second section we
introduce the functional framework. In the third and forth sections we prove
the existence theorem.

\section{Statement of the problem}

\subsection{The Klein-Gordon-Maxwell system}

Let $\psi:\mathbb{R}^{4}\rightarrow\mathbb{C}$ be a complex scalar field on
the space-time $\mathbb{R}^{4}$and $\Gamma$ be a 1- form on $\mathbb{R}^{4}$
whose coefficients $\Gamma_{j}$ are in the Lie algebra $u(1)$ of the group
$U(1)=S^{1},$ i.e. $\Gamma_{j}=-iA_{j}$, where $i$ is the imaginary unit and
$A_{j}$ ($j=0,..,3)$ are real maps defined in $\mathbb{R}^{4}.$

Consider the Abelian gauge theory related to $\psi$ and to $\Gamma$ and
described by the Lagrangian density (see e.g.\cite{yangL}, \cite{rub})
\begin{equation}
\mathcal{L}=\mathcal{L}_{0}+\mathcal{L}_{1}-W(\left|  \psi\right|
)\label{marisa}%
\end{equation}
where
\[
\mathcal{L}_{0}=-\frac{1}{2}\left\langle d_{A}\psi,d_{A}\psi\right\rangle
,\text{ }\mathcal{L}_{1}=-\frac{1}{2}\left\langle d_{A}A,d_{A}A\right\rangle
,\text{ }A=\overset{3}{\underset{j=0}{\sum}}A_{j}dx^{j}%
\]
and $W$ is a real $C^{1}$-function $.$ Here
\[
d_{A}=d-iA=\underset{}{\overset{3}{\underset{j=0}{\sum}}\left(  \frac
{\partial}{\partial x^{j}}-iA_{j}\right)  dx^{j}}%
\]
denotes the Weyl covariant differential and $\left\langle \cdot,\cdot
\right\rangle $ denotes the scalar product between forms with respect the
Minkowski metric in $\mathbb{R}^{4}$.

Since $A_{j}$ are real maps,
\[
d_{A}A=dA-iA\wedge A=dA
\]

Now we set
\[
\mathbf{%
A\mathbf{=(}%
}A_{1},A_{2},A_{3}\mathbf{%
)%
}\text{ and }\phi=-A_{0}%
\]
If we set $t=-x_{0}$ and $x=(x_{1},x_{2},x_{3})$ and use vector notation, the
Lagrangian densities $\mathcal{L}_{0},\mathcal{L}_{1}$ can be written as
follows
\begin{equation}
\mathcal{L}_{0}=\frac{1}{2}\left[  \left|  \left(  \partial_{t}+i\phi\right)
\psi\right|  ^{2}-\left|  \left(  \nabla-i\mathbf{A}\right)  \psi\right|
^{2}\right]  .
\end{equation}

\[
\mathcal{L}_{1}=\frac{1}{2}\left\vert \partial_{t}\mathbf{%
A%
}+\nabla\phi\right\vert ^{2}-\frac{1}{2}\left\vert \nabla\times\mathbf{A}%
\right\vert ^{2}%
\]
Here $\nabla\times$ and $\nabla$ denote respectively the curl and the gradient
operators with respect the $x$ variable and $\partial_{t}$ denotes the
derivative with respect to $t$ variable.

Observe that $\mathcal{L}_{1}$ is the Maxwell Lagrangian density of the
electromagnetic field
\[
\mathbf{E}=-\partial_{t}\mathbf{%
A%
}-\nabla\phi,\text{ }\mathbf{H=}\nabla\times\mathbf{A}%
\]
Now consider the total action of the Abelian gauge theory
\begin{equation}
\mathcal{S}=\int\left(  \mathcal{L}_{0}+\mathcal{L}_{1}-W(\left|  \psi\right|
)\right)  dxdt\label{ac1}%
\end{equation}

Making the variation of $\mathcal{S}$ with respect to $\psi,$ $\phi$ and
$\mathbf{A}$ we get the system of equations
\begin{equation}
\left(  \partial_{t}+i\phi\right)  ^{2}\psi-\left(  \nabla-i\mathbf{A}\right)
^{2}\psi+W^{^{\prime}}(\left|  \psi\right|  )\frac{\psi}{\left|  \psi\right|
}=0\label{e1}%
\end{equation}
\begin{equation}
\nabla\cdot\left(  \partial_{t}\mathbf{%
A%
}+\nabla\phi\right)  =\left(  \operatorname{Im}\frac{\partial_{t}\psi}{\psi
}+\phi\right)  \left|  \psi\right|  ^{2}\label{e2}%
\end{equation}
\begin{equation}
\nabla\times\left(  \nabla\times\mathbf{A}\right)  +\partial_{t}\left(
\partial_{t}\mathbf{%
A%
}+\nabla\phi\right)  =\left(  \operatorname{Im}\frac{\nabla\psi}{\psi
}-\mathbf{A}\right)  \left|  \psi\right|  ^{2}\label{e3}%
\end{equation}
Here $\nabla\cdot$ denotes the divergence operator. We recall that these
equations are gauge invariant i.e. if $\psi$, $\mathbf{A,}$ $\phi$ solve
(\ref{e1}), (\ref{e2}), (\ref{e3}), then also $e^{i\chi}\psi$, $\mathbf{A+}%
\nabla\chi\mathbf{,}$ $\phi-\partial_{t}\chi$ (with $\chi\in C^{\infty}\left(
\mathbb{R}^{4}\right)  $) solve the same equations.

\subsection{Stationary solutions and vortices}

We look for stationary solutions of (\ref{e1}), (\ref{e2}), (\ref{e3}), namely
solutions of the form
\begin{align*}
\psi\left(  t,x\right)   & =u\left(  x\right)  e^{i\left(  S(x)-\omega
t\right)  },\text{\ }u\in\mathbb{R}^{+},\ \omega\in\mathbb{R},\text{\ }%
S\in\frac{\mathbb{R}}{2\pi\mathbb{Z}}\\
\partial_{t}\mathbf{A}  & =0\mathbf{,\ }\partial_{t}\phi=0
\end{align*}
A stationary solution solves the following set of equations:
\begin{equation}
-\Delta u+\left[  \left|  \nabla S-\mathbf{A}\right|  ^{2}-\left(  \phi
-\omega\right)  ^{2}\right]  \,u+W^{\prime}\left(  u\right)  =0\label{h1}%
\end{equation}
\begin{equation}
-\nabla\cdot\left[  \left(  \nabla S-\mathbf{A}\right)  u^{2}\right]
=0\label{h2}%
\end{equation}
\begin{equation}
-\Delta\phi=\left(  \omega-\phi\right)  u^{2}\;\label{h3}%
\end{equation}
\begin{equation}
\nabla\times\left(  \nabla\times\mathbf{A}\right)  =\left(  \nabla
S-\mathbf{A}\right)  u^{2}\;\label{h4}%
\end{equation}
Observe that equation (\ref{h2}) easily follows from equation (\ref{h4}). Then
we are reduced to study the system (\ref{h1}), (\ref{h3}), (\ref{h4}). Clearly
when $u=0,$ the only finite energy gauge potentials which solve (\ref{h3}),
(\ref{h4}) are the trivial ones $\mathbf{A=0},$ $\phi=0.$

It is possible to have three types of stationary non trivial solutions:

\begin{itemize}
\item electrostatic solutions: $\mathbf{A}=0$, $\phi\neq0;$

\item magneto-static solutions: $\mathbf{A}\neq0$, $\phi=0;$

\item electro-magneto-static solutions: $\mathbf{A}\neq0$, $\phi\neq0$.
\end{itemize}

Under suitable assumptions, all these types of solutions exist. The existence
and the non existence of electrostatic solutions for a system like
(\ref{h1}),..., (\ref{h3}) has been proved under different assumptions on $W$
(see \cite{bf} , \cite{befo}, \cite{ca}, \cite{tea}, \cite{tea2}, \cite{dav}).
In particular the existence of radially symmetric, finite energy electrostatic
solutions has been analyzed.

Here we are interested in magneto-static and electro-magneto-static solutions,
in particular we shall study the existence of vortices in the sense of the
definition stated below. We set
\[
\Sigma=\left\{  \left(  x_{1},x_{2},x_{3}\right)  \in\mathbb{R}^{3}%
:x_{1}=x_{2}=0\right\}
\]
and we define the map
\[
\theta:\mathbb{R}^{3}\backslash\Sigma\rightarrow\frac{\mathbb{R}}%
{2\pi\mathbb{Z}}%
\]
\[
\theta(x_{1},x_{2},x_{3})=\operatorname{Im}\log(x_{1}+ix_{2}).
\]
Since there are three equations (\ref{h1}), (\ref{h3}), (\ref{h4}) in the
unknowns $u,$ $S$, $\mathbf{A}$, $\phi,$ we shall take $S=k\vartheta$ ($k$
integer) and we shall solve with respect $u,$ $\mathbf{A}$, $\phi.$ So we give
the following definition. A solution of Eq. (\ref{h1}), (\ref{h3}), (\ref{h4})
is called vortex if $\psi$ has the following form
\begin{equation}
\psi(t,x)=u(x)\,e^{i\left(  k\theta(x)-\omega t\right)  };\ k\in
\mathbb{Z-}\left\{  0\right\}  .\label{ans}%
\end{equation}
We shall see in Theorem \ref{main} that if $\psi$ has the above form, then
$\mathbf{A}$ and $\mathbf{H}$ look like the fields created by a finite solenoid.

Observe that $\theta\in C^{\infty}\left(  \mathbb{R}^{3}\backslash\Sigma
,\frac{\mathbb{R}}{2\pi\mathbb{Z}}\right)  $ and$\ \nabla\theta\in C^{\infty
}\left(  \mathbb{R}^{3}\backslash\Sigma,\mathbb{R}^{3}\right)  ,$ namely
\[
\nabla\theta(x)=\left(  \frac{x_{2}}{x_{1}^{2}+x_{2}^{2}},\ \frac{-x_{1}%
}{x_{1}^{2}+x_{2}^{2}},\ 0\right)  .
\]

Using this ansatz equations (\ref{h1}), (\ref{h3}), (\ref{h4}) become
\begin{equation}
-\Delta u+\left[  \left|  k\nabla\theta-\mathbf{A}\right|  ^{2}-\left(
\phi-\omega\right)  ^{2}\right]  \,u+W^{\prime}(u)=0\label{z1}%
\end{equation}
\begin{equation}
-\Delta\phi=\left(  \omega-\phi\right)  u^{2}\;\label{z3}%
\end{equation}
\begin{equation}
\nabla\times\left(  \nabla\times\mathbf{A}\right)  =\left(  k\nabla
\theta-\mathbf{A}\right)  u^{2}\;\label{z4}%
\end{equation}

It can be shown (see \cite{befo}) that the energy of a solution $(u,\phi
,\mathbf{A)}$ of equations (\ref{z1}), (\ref{z3}), (\ref{z4}) has the
following expression
\begin{equation}
\mathcal{E}=\frac{1}{2}\int\left(  \left|  \nabla u\right|  ^{2}+\left|
\nabla\phi\right|  ^{2}+\left|  \nabla\mathbf{A}\right|  ^{2}+(\left|
\mathbf{A}-k\nabla\theta\right|  ^{2}+\left(  \phi-\omega\right)  ^{2}%
)\,u^{2}\right)  +\int W(u)\nonumber
\end{equation}

\subsection{The main results}

Let
\begin{equation}
W(s)=\frac{1}{2}s^{2}-F(s)\label{n}%
\end{equation}
where $F$ is a $C^{2}$ real function satisfying the following assumptions:
\begin{equation}
F(0)=F^{\prime}(0)=F^{\prime\prime}(0)=0\label{nn}%
\end{equation}

There are constants $c>0$ and $p$ with $2<p<6$ such that
\begin{equation}
\left\vert F^{\prime}(s)\right\vert \leq cs^{p-1},\text{ }s\geq0\label{nnn}%
\end{equation}
\ \ \ \ \ \ \ \ \ \
\begin{equation}
sF^{\prime}(s)\geq pF(s)>0\text{ for }s>0\label{nnnn}%
\end{equation}
Thus a typical function $W$ satisfying our assumptions is
\begin{equation}
W(s)=\frac{s^{2}}{2}-\frac{s^{p}}{p}\text{ }s\geq0,\text{ }2<p<6.\label{tip}%
\end{equation}

Moreover, for technical reasons it is useful to assume that $W$ is defined for
all $s\in\mathbb{R}$ just setting
\[
W(s)=\frac{s^{2}}{2}\ for\ s<0
\]

We shall prove the following existence result for vortex solutions

\begin{theorem}
\label{main}Assume that the function $W$ satisfy assumptions (\ref{n}%
)...(\ref{nnnn}) with $2<p<6$ and set
\[
\omega_{p}=\min\left( 1,\sqrt{\frac{p-2}{2}}\right) \text{ }%
\]
Then for any $\omega\in\left(  -\omega_{p},\omega_{p}\right)  \ $and any
$k\in\mathbb{Z}$ \ the equations (\ref{z1}), (\ref{z3}), (\ref{z4}) admit a
solution (in the sense of distributions on $\mathbb{R}^{3}$) $(u,\phi
\mathbf{,A)}$ with $u\neq0$ which satisfy the following properties:

\begin{itemize}
\item (a) $\int\left[  \left|  \nabla u\right|  ^{2}+\left(  1+\frac{1}{r^{2}%
}\right)  u^{2}\right]  dx<+\infty$ , $r^{2}=x_{1}^{2}+x_{2}^{2}$

\item (b) $\int\left|  \nabla\phi\right|  ^{2}+\left|  \nabla\mathbf{A}%
\right|  ^{2}<+\infty$

\item (c) $u\geq0$

\item (d) there exists a real function $b$ such that $\mathbf{A}=b\nabla
\theta$.

\item (e) $u,$ $\phi$ and $\left\vert \mathbf{A}\right\vert $ have cylindrical
symmetry, i.e. they depend only on $r$ and $x_{3}$
\end{itemize}

\noindent Moreover

i) if $\omega\neq0$ and $k=0,$ then $\phi\neq0$ and $\mathbf{A}$ $=0$
\textit{(electrostatic solutions).}

ii) if $\omega=0$ and $k\neq0,$ then $\phi=0$ and $\mathbf{A}$ $\neq0$
\textit{(magnetostatic vortices).}

iii) if $\omega\neq0$ and $k\neq0,$ then $\phi\neq0$ and $\mathbf{A}$ $\neq0$
\textit{(electromagnetostatic vortices)}
\end{theorem}

\begin{remark}
The properties (a) and (b) guarantee that the energy is finite.
\end{remark}

\begin{remark}
\label{iv} If $\left(  u,\phi\mathbf{,A}\right)  $ with $u\neq0$ solves
(\ref{z1}), (\ref{z3}), (\ref{z4}), then assertions i), ii), iii) in Theorem
\ref{main} follow immediately from (\ref{z3}), (\ref{z4}).
\end{remark}

\begin{remark}
The \textit{magnetostatic vortices are the critical points of the functional
}
\[
\int\frac{1}{2}\left|  \left(  \nabla-i\mathbf{A}\right)  \psi\right|
^{2}+\frac{1}{2}\left|  \nabla\times\mathbf{A}\right|  ^{2}+W(\left|
\psi\right|  )
\]
Thus, they can be interpreted as solutions of a Euclidean gauge theory in
dimension $3.$
\end{remark}

\begin{remark}
\label{v} By the presence of the term $\nabla\theta$ equations (\ref{z1}),
(\ref{z4}) are not invariant under the $O(3)$ group action as it happens for
the equations (\ref{e1}), (\ref{e2}), (\ref{e3}) we started from. Indeed there
is a breaking of radial symmetry and the solutions $u$, $\phi$ and $\mathbf{A}
$ have only a $O(1)=S^{1}$ symmetry (see (d) and (e) in Theorem \ref{main}).
\end{remark}

Finally we point out that the solutions in Theorem \ref{main} corresponding to
different $k$ or different $\omega$ cannot be obtained one from the other by
means of a smooth gauge transformation. In fact the following proposition holds

\begin{proposition}
\label{inv}Let $\left(  u_{1}(x)\,e^{i\left(  k_{1}\theta(x)-\omega
_{1}t\right)  },\phi_{1},\mathbf{A}_{1}\right)  $ and $\left(  u_{2}%
(x)\,e^{i\left(  k_{2}\theta(x)-\omega_{2}t\right)  },\phi_{2},\mathbf{A}%
_{2}\right)  $ be two vortex solutions with $\left(  \phi_{1},\mathbf{A}%
_{1}\right)  ,\left(  \phi_{2},\mathbf{A}_{2}\right)  \in\mathcal{D}%
^{1,2}\times\left(  \mathcal{D}^{1,2}\right)  ^{3}$ (see (\ref{inner}))$.$
Assume that $\omega_{2}\neq\omega_{1}$ or $k_{2}\neq k_{1}.$ Then these
solutions cannot be obtained one from the other by means of a gauge
transformation.$.$
\end{proposition}

\begin{proof}
By assumption there exists a map $\chi$ such that
\begin{equation}
u_{1}(x)\,e^{i\left(  k_{1}\theta(x)-\omega_{1}t+\chi\right)  }=u_{2}%
(x)\,e^{i\left(  k_{2}\theta(x)-\omega_{2}t\right)  }\label{la}%
\end{equation}

\begin{equation}
\phi_{1}-\partial_{t}\chi=\phi_{2},\text{ }\mathbf{A}_{1}+\nabla
\chi=\mathbf{A}_{2}\label{le}%
\end{equation}

then by (\ref{la}) we have
\begin{equation}
u_{1}=u_{2},\text{ }\chi(x,t)=\left(  k_{2}-k_{1}\right)  \theta(x)-\left(
\omega_{2}-\omega_{1}\right)  t\label{li}%
\end{equation}

from which
\begin{equation}
\partial_{t}\chi=\omega_{1}-\omega_{2},\text{ }\nabla\chi=\left(  k_{2}%
-k_{1}\right)  \nabla\theta(x)\label{lo}%
\end{equation}
By (\ref{lo}), (\ref{le}) we get
\begin{align}
\phi_{2}-\phi_{1}  & =\omega_{2}-\omega_{1}\label{ca}\\
\mathbf{A}_{2}-\mathbf{A}_{1}  & =\left(  k_{2}-k_{1}\right)  \nabla
\theta(x)\label{ce}%
\end{align}
Since $\left(  \phi_{1},\mathbf{A}_{1}\right)  ,\left(  \phi_{2}%
,\mathbf{A}_{2}\right)  \in\mathcal{D}^{1,2}\times\left(  \mathcal{D}%
^{1,2}\right)  ^{3}, $ by (\ref{ca}) and (\ref{ce}) we deduce $\omega
_{2}=\omega_{1}$ and $k_{2}=k_{1}$ and this contradicts our assumptions.
\end{proof}

\section{The functional setting}

\subsection{Weak solutions}

Let $H^{1}$ denote the usual Sobolev space with norm
\[
\left\Vert u\right\Vert _{H^{1}}^{2}=\int(\left\vert \nabla u\right\vert
^{2}+u^{2})dx;
\]
moreover we need to use also the weighted Sobolev space $\hat{H}^{1}$ whose
norm is given by
\[
\left\Vert u\right\Vert _{\hat{H}^{1}}^{2}=\int\left[  \left\vert \nabla
u\right\vert ^{2}+\left(  1+\frac{1}{r^{2}}\right)  u^{2}\right]  dx
\]
where $r=\sqrt{x_{1}^{2}+x_{2}^{2}}.$

We set $\mathcal{D}=C_{0}^{\infty}(\mathbb{R}^{3})$ and we denote by
$\mathcal{D}^{1,2}$ the completion of $\mathcal{D}$ with respect to the inner
product
\begin{equation}
\left(  v\mid w\right)  _{\mathcal{D}^{1,2}}=\int\nabla v\cdot\nabla
wdx\label{inner}%
\end{equation}
Here and in the following the dot $\cdot$ will denote the Euclidean inner
product in $\mathbb{R}^{3}.$

We set
\[
H=\hat{H}^{1}\times\mathcal{D}^{1,2}\times\left(  \mathcal{D}^{1,2}\right)
^{3}%
\]
\begin{equation}
\left\Vert \left(  u,\phi,\mathbf{A}\right)  \right\Vert _{H}^{2}%
=\int\left\vert \nabla u\right\vert ^{2}+\left(  1+\frac{1}{r^{2}}\right)
u^{2}+\left\vert \nabla\phi\right\vert ^{2}+\left\vert \nabla\mathbf{A}%
\right\vert ^{2}.\label{usual}%
\end{equation}

Now we consider the functional
\begin{align}
J(u,\phi\mathbf{,A})  & =\frac{1}{2}\int\left|  \nabla u\right|  ^{2}-\left|
\nabla\phi\right|  ^{2}+\left|  \nabla\times\mathbf{A}\right|  ^{2}\nonumber\\
& +\frac{1}{2}\int\left[  \left|  \mathbf{A}-k\nabla\theta\right|
^{2}-\left(  \phi-\omega\right)  ^{2}\right]  \,u^{2}+\int
W(u)\label{functional}%
\end{align}

\noindent where $\left(  u,\phi,\mathbf{A}\right)  \in H.$ The equations
(\ref{z1}),\ (\ref{z3}) and (\ref{z4}) are the Euler-Lagrange equations of the
functional $J$. In fact the following lemma holds:

\ 

\begin{lemma}
\label{tecnico}The functional $J$ is $C^{1}$ on $H$.
\end{lemma}

$.$

\begin{proof}
The variations of $J$ with respect to $u,\phi,\mathbf{A}$ are respectively
\begin{equation}
J_{u}^{\prime}=-\Delta u+\left[  \left\vert k\nabla\theta-\mathbf{A}%
\right\vert ^{2}-\left(  \phi-\omega\right)  ^{2}\right]  \,u+W^{\prime
}\left(  u\right) \label{be1}%
\end{equation}%
\begin{equation}
J_{\phi}^{\prime}=\Delta\phi-\left(  \phi-\omega\right)  u^{2}\;\label{be2}%
\end{equation}%
\begin{equation}
J_{\mathbf{A}}^{\prime}=\nabla\times\nabla\times\mathbf{A}+\left(
\mathbf{A}-k\nabla\theta\right)  u^{2}\;\label{be3}%
\end{equation}
We need to prove that the maps
\begin{align}
\left(  u,\phi,\mathbf{A}\right)   & \in H\mapsto J_{u}^{\prime}\left(
u,\phi,\mathbf{A}\right)  \in\left(  \hat{H}^{1}\right)  ^{\prime}\label{n1}\\
\left(  u,\phi,\mathbf{A}\right)   & \in H\mapsto J_{\phi}^{\prime}\left(
u,\phi,\mathbf{A}\right)  \in(\mathcal{D}^{1,2})^{\prime}\label{n2}\\
\left(  u,\phi,\mathbf{A}\right)   & \in H\mapsto J_{\mathbf{A}}^{\prime
}\left(  u,\phi,\mathbf{A}\right)  \in\left(  (\mathcal{D}^{1,2})^{3}\right)
^{\prime}\label{n3}%
\end{align}
are continuous. First we show that the map
\begin{equation}
\left(  u,\mathbf{A}\right)  \rightarrow\left\vert \mathbf{A}\right\vert
^{2}u\label{tec1}%
\end{equation}
is continuous from $\hat{H}^{1}\times(\mathcal{D}^{1,2})^{3}$ to the Lebesgue
space $L^{\frac{6}{5}\text{ }}$(and therefore to $\left(  \hat{H}^{1}\right)
^{\prime}$).

In fact, if $\left(  u_{n},\mathbf{A}_{n}\right)  $ $\rightarrow\left(
u,\mathbf{A}\right)  $ in $\hat{H}^{1}\times(\mathcal{D}^{1,2})^{3}$, then, by
Sobolev embedding and H\"{o}lder inequalities, we easily obtain:
\[
\left\|  \left|  \mathbf{A}_{n}\right|  ^{2}u_{n}-\left|  \mathbf{A}\right|
^{2}u\right\|  _{L^{\frac{6}{5}}}\leq
\]
\[
\left\|  \left|  \mathbf{A}_{n}\right|  ^{2}(u_{n}-u)\right\|  _{L^{\frac
{6}{5}}}+\left\|  \left(  \left|  \mathbf{A}_{n}\right|  ^{2}-\left|
\mathbf{A}\right|  ^{2}\right)  u\right\|  _{L^{\frac{6}{5}}}\leq
\]
\[
\leq\left\|  \mathbf{A}_{n}\right\|  _{L^{6}}^{2}\left\|  u_{n}-u\right\|
_{L^{2}}+\left\|  \left|  \mathbf{A}_{n}\right|  ^{2}-\left|  \mathbf{A}%
\right|  ^{2}\right\|  _{L^{3}}\left\|  u\right\|  _{L^{2}}%
\]
from which we get the conclusion.

Analogously it can be shown that the maps
\begin{equation}
\left(  u,\phi\right)  \in\hat{H}^{1}\times\mathcal{D}^{1,2}\rightarrow
\phi^{2}u\in L^{\frac{6}{5}\text{ }}\subset\left(  \hat{H}^{1}\right)
^{\prime}\label{tec2}%
\end{equation}

$,$
\begin{equation}
\left(  u,\mathbf{A}\right)  \in\hat{H}^{1}\times(\mathcal{D}^{1,2}%
)^{3}\rightarrow\mathbf{A}u^{2}\in(L^{\frac{6}{5}\text{ }})^{3}\subset\left(
(\mathcal{D}^{1,2})^{3}\right)  ^{\prime}\label{tec3}%
\end{equation}

\begin{equation}
\left(  u,\phi\right)  \in\hat{H}^{1}\times\mathcal{D}^{1,2}\rightarrow\phi
u^{2}\in L^{\frac{6}{5}\text{ }}\subset(\mathcal{D}^{1,2})^{\prime
}\label{tec4}%
\end{equation}
are continuous. Moreover it is immediate to see that the linear map
\begin{equation}
u\in\hat{H}^{1}\rightarrow\left|  \nabla\theta\right|  ^{2}u=\frac{u}{r^{2}%
}\in\left(  \hat{H}^{1}\right)  ^{\prime}\label{tec5}%
\end{equation}
is continuous.

Also the map
\begin{equation}
\left(  u,\mathbf{A}\right)  \in\hat{H}^{1}\times\left(  \mathcal{D}%
^{1,2}\right)  ^{3}\rightarrow u\mathbf{A\cdot}\nabla\theta\in L^{\frac{3}{2}%
}\subset\left(  \hat{H}^{1}\right)  ^{\prime}\label{tec6}%
\end{equation}
is continuous$.$ In fact, let $\left(  u_{n},\mathbf{A}_{n}\right)
\rightarrow\left(  u,\mathbf{A}\right)  $ in $\hat{H}^{1}\times\left(
\mathcal{D}^{1,2}\right)  ^{3}.$ For all $\varphi\in H^{1},$ we have that
\[
\left|  \int(u_{n}\mathbf{A}_{n}\mathbf{\cdot}\nabla\theta-u\mathbf{A\cdot
}\nabla\theta)\varphi\right|  \leq
\]
\[
\leq\left|  \int u_{n}\left(  \mathbf{A}_{n}-\mathbf{A}\right)  \mathbf{\cdot
}\nabla\theta\varphi\right|  +\left|  \int\mathbf{A\cdot}\nabla\theta
(u_{n}-u)\varphi\right|  \leq
\]
\[
\leq\int\frac{\left|  u_{n}\varphi\right|  \,\left|  \mathbf{A}_{n}%
-\mathbf{A}\right|  }{r}+\int\frac{\left|  \left(  u_{n}-u\right)
\varphi\right|  \,\left|  \mathbf{A}\right|  }{r}\leq
\]
\[
\leq\left\|  u_{n}\right\|  _{\hat{H}^{1}}\left\|  \varphi\right\|  _{L^{3}%
}\left\|  \mathbf{A}_{n}-\mathbf{A}\right\|  _{L^{6}}+\left\|  u_{n}%
-u\right\|  _{\hat{H}^{1}}\left\|  \varphi\right\|  _{L^{3}}\left\|
\mathbf{A}\right\|  _{L^{6}}%
\]
from which the conclusion easily follows.

Arguing as before it is possible to show that also the map
\begin{equation}
u\in\hat{H}^{1}\rightarrow u^{2}\nabla\theta\in\left(  (\mathcal{D}^{1,2}%
)^{3}\right)  ^{\prime}\label{tec7}%
\end{equation}
is continuous. So by the continuity of (\ref{tec1}),...,(\ref{tec7}) we deduce
that (\ref{n1}), (\ref{n2}) and (\ref{n3}) are continuous.
\end{proof}

\bigskip

By the above lemma it follows that the critical points $\left(  u,\phi
,\mathbf{A}\right)  \in H$ of $J$ are weak solutions of eq. (\ref{z1}%
),\ (\ref{z3}) and (\ref{z4}), namely
\begin{equation}
\int\nabla u\cdot\nabla v+\left[  \left|  \mathbf{A}-k\nabla\theta\right|
^{2}-\left(  \phi-\omega\right)  ^{2}\right]  \,uv+W^{\prime}\left(  u\right)
v=0,\;\forall v\in\hat{H}^{1}\label{w1}%
\end{equation}
\begin{equation}
\int\nabla\phi\cdot\nabla w+\left(  \phi-\omega\right)  u^{2}w=0,\;\forall
w\in\mathcal{D}^{1,2}\label{we}%
\end{equation}
\begin{equation}
\int\nabla\mathbf{A}\cdot\nabla\mathbf{V}+u^{2}\left(  \mathbf{A}%
-k\nabla\theta\right)  \cdot\mathbf{V}=0,\;\forall\mathbf{V}\in(\mathcal{D}%
^{1,2})^{3}.\label{w3}%
\end{equation}

\subsection{Solutions in the sense of distributions}

\bigskip

Since $\mathcal{D}$ is not contained in $\hat{H}^{1},$ a solution $\left(
u,\phi,\mathbf{A}\right)  \in H$ of \ (\ref{w1}), (\ref{we}), (\ref{w3})) need
not to be a solution of (\ref{z1}), (\ref{z3}), (\ref{z4}) in the sense of
distributions on $\mathbb{R}^{3}$. In fact, since $\nabla\theta\left(
x\right)  $ is singular on $\Sigma,$ it might be that for some test function
$v\in\mathcal{D},$ the integral $\int\left\vert \mathbf{A}-k\nabla
\theta\right\vert ^{2}\,uv$ diverges, unless $u$ is sufficiently small as
$x\rightarrow\Sigma.$

In this section we will show that this fact does not occur, namely the
singularity is removable in the sense of the following theorem:

\begin{theorem}
\label{finale}Let $(u_{0},\phi_{0}\mathbf{,A}_{0}\mathbf{)}$ $\in H$ be a
solution of (\ref{w1}), (\ref{we}), (\ref{w3}) (i.e. a critical point of $J).
$ Then $(u_{0},\phi_{0}\mathbf{,A}_{0}\mathbf{)}$ is a solution of eq.
(\ref{z1}),\ (\ref{z3}) and (\ref{z4}) in the sense of distribution, namely
\begin{equation}
\int\nabla u_{0}\cdot\nabla v+\left[  \left\vert \mathbf{A}_{0}-k\nabla
\theta\right\vert ^{2}-\left(  \phi_{0}-\omega\right)  ^{2}\right]
\,u_{0}v+W^{\prime}\left(  u_{0}\right)  v=0,\;\forall v\in\mathcal{D}%
\label{im}%
\end{equation}%
\begin{equation}
\int\nabla\phi_{0}\cdot\nabla w+\left(  \phi_{0}-\omega\right)  u_{0}%
^{2}w=0,\;\forall w\in\mathcal{D}\label{dis}%
\end{equation}%
\begin{equation}
\int\nabla\mathbf{A}_{0}\cdot\nabla\mathbf{V}+\left(  \mathbf{A}_{0}%
-k\nabla\theta\right)  u_{0}^{2}\cdot\mathbf{V}=0,\;\forall\mathbf{V}%
\in\mathcal{D}^{3}.\label{um}%
\end{equation}

\end{theorem}

\bigskip Before proving Theorem \ref{finale} we need some lemmas.

\begin{lemma}
\label{pis}Let $\left(  u_{0},\phi_{0},\ \mathbf{A}\right)  \ $be a solution
of (\ref{w1}), (\ref{we}), (\ref{w3}). Then $u_{0}\left(  x\right)  \geq0$ a.e.
\end{lemma}

\begin{proof}
Set $u_{0}^{-}(x)=\min\left\{  u_{0}(x),0\right\}  .$ Then, if in eq.
(\ref{w1}) we set $v=u_{0}^{-},$ we get
\begin{equation}
\int\left|  \nabla u_{0}^{-}\right|  ^{2}+\left[  \left|  \mathbf{A}%
-k\nabla\theta\right|  ^{2}-\left(  \phi_{0}-\omega\right)  ^{2}\right]
\,\left(  u_{0}^{-}\right)  ^{2}+W^{\prime}(u_{0}^{-})u_{0}^{-}=0\label{g}%
\end{equation}

and since $W^{\prime}(s)=s$ for $s\leq0,$ we get
\begin{equation}
\int\left\vert \nabla u_{0}^{-}\right\vert ^{2}+\left[  \left\vert
\mathbf{A}-k\nabla\theta\right\vert ^{2}+1-\left(  \phi_{0}-\omega\right)
^{2}\right]  \,\left(  u_{0}^{-}\right)  ^{2}=0\label{h}%
\end{equation}

We shall show (see Lemma \ref{massimo}) that
\[
0\leq\frac{\phi_{0}}{\omega}\leq1
\]
So, since $\omega^{2}<1,$ we have $\left(  \phi_{0}-\omega\right)  ^{2}\leq1.$
Then
\begin{equation}
\int\left[  1-\left(  \phi_{0}-\omega\right)  ^{2}\right]  \left(  u_{0}%
^{-}\right)  ^{2}\geq0\label{i}%
\end{equation}

From (\ref{h}) and (\ref{i}) we deduce that $u_{0}^{-}=0.$
\end{proof}

Let $\chi_{n}$ ($n$ positive integer) be a family of smooth functions
depending only on $r=\sqrt{x_{1}^{2}+x_{2}^{2}}$ and $x_{3}$ and which satisfy
the following assumptions:

\begin{itemize}
\item $\chi_{n}\left(  r,x_{3}\right)  =1$ for $r\geq\frac{2}{n}$

\item $\chi_{n}\left(  r,x_{3}\right)  =0$ for $r\leq\frac{1}{n}$

\item $\left|  \chi_{n}\left(  r,x_{3}\right)  \right|  \leq1$

\item $\left|  \nabla\chi_{n}\left(  r,x_{3}\right)  \right|  \leq2n$

\item $\chi_{n+1}\left(  r,x_{3}\right)  \geq\chi_{n}\left(  r,x_{3}\right)  $
\end{itemize}

\begin{lemma}
\label{caccona}Let $\varphi$ be a function in $H^{1}\cap L^{\infty}$ with
bounded support and set $\varphi_{n}=\varphi\cdot\chi_{n}.$ Then, up to a
subsequence, we have that
\[
\varphi_{n}\rightarrow\varphi\text{ weakly in }H^{1}%
\]

\end{lemma}

\begin{proof}
Clearly $\varphi_{n}\rightarrow\varphi$ $a.e.$ Then, by standard arguments,
the conclusion holds if we show that $\left\{  \varphi_{n}\right\}  $ is
bounded in $H^{1}.$ Clearly $\left\{  \varphi_{n}\right\}  $ is bounded in
$L^{2}.$ Let us now prove that
\[
\left\{  \int\left|  \nabla\varphi_{n}\right|  ^{2}\right\}  \text{ is
bounded}%
\]
We have
\begin{align*}
\int\left|  \nabla\varphi_{n}\right|  ^{2}  & \leq2\int\left|  \nabla
\varphi\cdot\chi_{n}\right|  ^{2}+\left|  \varphi\cdot\nabla\chi_{n}\right|
^{2}\\
& \leq2\int\left|  \nabla\varphi\right|  ^{2}+2\int_{\Gamma_{\varepsilon}%
}\left|  \varphi\cdot\nabla\chi_{n}\right|  ^{2}%
\end{align*}
where
\[
\Gamma_{\varepsilon}=\left\{  x\in\mathbb{R}^{3}:\varphi\neq0\ and\ \left|
\nabla\chi_{n}\left(  r,z\right)  \right|  \neq0\right\}
\]
By our construction, $\left|  \Gamma_{\varepsilon}\right|  \leq c/n^{2}$ where
$c$ depends only on $\varphi.$ Thus
\begin{align*}
\int\left|  \nabla\varphi_{n}\right|  ^{2}  & \leq2\int\left|  \nabla
\varphi\right|  ^{2}+2\left\|  \varphi\right\|  _{L^{\infty}}^{2}\int
_{\Gamma_{\varepsilon}}\left|  \nabla\chi_{n}\right|  ^{2}\\
& \leq2\int\left|  \nabla\varphi\right|  ^{2}+2\left\|  \varphi\right\|
_{L^{\infty}}^{2}\cdot\left|  \Gamma_{\varepsilon}\right|  \cdot\left\|
\nabla\chi_{n}\right\|  _{L^{\infty}}^{2}\\
& \leq2\int\left|  \nabla\varphi\right|  ^{2}+8c\left\|  \varphi\right\|
_{L^{\infty}}^{2}%
\end{align*}
Thus, $\varphi_{n}$ is bounded in $H^{1}$ and $\varphi_{n}\rightarrow\varphi$
weakly in $H^{1}.$

\bigskip
\end{proof}

\bigskip

\bigskip Now we are ready to prove Theorem \ref{finale}

\begin{proof}
Clearly (\ref{dis}) and (\ref{um}) immediately follows by (\ref{we}) and
(\ref{w3}). Let us prove (\ref{im}). We take any $v\in\mathcal{D}$ and set
$\varphi_{n}=v^{+}\chi_{n}$. Then taking $\varphi_{n}$ as test function in Eq.
(\ref{w1}), we have
\begin{equation}
\int\nabla u_{0}\cdot\nabla\varphi_{n}+\left[  \left|  \mathbf{A}_{0}%
-k\nabla\theta\right|  ^{2}-\left(  \phi_{0}-\omega\right)  ^{2}\right]
\,u_{0}\varphi_{n}+W^{\prime}\left(  u_{0}\right)  \varphi_{n}=0\ \label{ecco}%
\end{equation}

Equation (\ref{ecco}) can be written as follows
\begin{equation}
A_{n}+B_{n}+C_{n}+D_{n}=0\label{simbol}%
\end{equation}

where
\begin{equation}
A_{n}=\int\nabla u_{0}\cdot\nabla\varphi_{n},\text{ }B_{n}=\int\left(
\mathbf{A}_{0}^{2}u_{0}-\left(  \phi_{0}-\omega\right)  ^{2}u_{0}+W^{\prime
}\left(  u_{0}\right)  \right)  \varphi_{n}\label{def1}%
\end{equation}

\begin{equation}
C_{n}=-2\int\mathbf{A}_{0}\cdot k\nabla\theta\ u_{0}\varphi_{n},\text{ }%
D_{n}=\int\left|  k\nabla\theta\right|  ^{2}u_{0}\varphi_{n}\label{def2}%
\end{equation}
By Lemma \ref{caccona}
\begin{equation}
\varphi_{n}\rightarrow v^{+}\text{ weakly in }H^{1}\label{w}%
\end{equation}

Then we have
\begin{equation}
A_{n}\rightarrow\int\nabla u_{0}\cdot\nabla v^{+}\label{pa}%
\end{equation}

Now
\[
\left(  \mathbf{A}_{0}^{2}u_{0}-\left(  \phi_{0}-\omega\right)  ^{2}%
u_{0}+W^{\prime}\left(  u_{0}\right)  \right)  \in L^{6/5}=\left(
L^{6}\right)  ^{\prime}%
\]
Then, using again (\ref{w}) and by the embedding $H^{1}\subset L^{6}$, we
have
\begin{equation}
B_{n}\rightarrow\int\left(  \mathbf{A}_{0}^{2}u_{0}-\left(  \phi_{0}%
-\omega\right)  ^{2}u_{0}+W^{\prime}\left(  u_{0}\right)  \right)
v^{+}<\infty\label{bi}%
\end{equation}

Now we shall prove that
\begin{equation}
C_{n}\rightarrow-2\int\mathbf{A}_{0}\cdot k\nabla\theta\ u_{0}v^{+}%
<\infty\label{ba}%
\end{equation}
Set
\[
C=B_{R}\times\left[  -l,l\right]  ,\text{ }B_{R}=\left\{  (x_{1},x_{2}%
)\in\mathbb{R}^{2}:r^{2}=x_{1}^{2}+x_{2}^{2}<R\right\}
\]
where $l,$ $R>0$ are so large that the cylinder $C$ contains the support of
$v^{+}.$

Then
\begin{align}
\int\left(  \frac{\varphi_{n}}{r}\right)  ^{\frac{3}{2}}dx  & =\int_{C}\left(
\frac{v^{+}\chi_{n}}{r}\right)  ^{\frac{3}{2}}dx\\
& \leq c_{1}\int_{-l}^{l}\int_{0}^{R}\left(  \frac{1}{r}\right)  ^{\frac{3}%
{2}}rdrdx_{3}=M<\infty\label{pos}%
\end{align}

where $c_{1}=2\pi\sup\left(  v^{+}\right)  ^{\frac{3}{2}}.$ By (\ref{pos}) we
have
\begin{equation}
\int\left|  \mathbf{A}_{0}\cdot\nabla\theta\ u_{0}\varphi_{n}\right|
dx\leq\left\|  u_{0}\mathbf{A}_{0}\ \right\|  _{L^{3}}\left\|  \frac
{\varphi_{n}}{r}\right\|  _{L^{\frac{3}{2}}}\leq\left\|  u_{0}\mathbf{A}%
_{0}\ \right\|  _{L^{3}}M^{\frac{2}{3}}\label{posse}%
\end{equation}

Now
\[
\left\vert \mathbf{A}_{0}\cdot\nabla\theta\ u_{0}\varphi_{n}\right\vert
\rightarrow\left\vert \mathbf{A}_{0}\cdot\nabla\theta\ u_{0}v^{+}\right\vert
\text{ a.e. in }\mathbb{R}^{3}%
\]
and the sequence $\left\{  \left\vert \mathbf{A}_{0}\cdot\nabla\theta
\ u_{0}\varphi_{n}\right\vert \right\}  $ is monotone $.$ Then, by the
monotone convergence theorem, we get
\begin{equation}
\int\left\vert \mathbf{A}_{0}\cdot k\nabla\theta\ u_{0}\varphi_{n}\right\vert
dx\rightarrow\int\left\vert \mathbf{A}_{0}\cdot k\nabla\theta\ u_{0}%
v^{+}\right\vert dx\label{from}%
\end{equation}

By (\ref{posse}) and (\ref{from}) we deduce that
\begin{equation}
\int\left|  \mathbf{A}_{0}\cdot k\nabla\theta\ u_{0}v^{+}\right|
dx<\infty\label{byy}%
\end{equation}

Then, since
\[
\left|  \mathbf{A}_{0}\cdot\nabla\theta\ u_{0}\varphi_{n}\right|  \leq\left|
\mathbf{A}_{0}\cdot\nabla\theta\ u_{0}v^{+}\right|  \in L^{1}%
\]
by the dominated convergence Theorem, we get (\ref{ba}). Finally we prove
that
\begin{equation}
D_{n}\rightarrow\int\left|  k\nabla\theta\right|  ^{2}u_{0}v^{+}%
<\infty\label{fin}%
\end{equation}
By (\ref{simbol}), (\ref{pa}), (\ref{bi}) and (\ref{ba}) we have that
\begin{equation}
D_{n}=\int\left|  k\nabla\theta\right|  ^{2}u_{0}\varphi_{n}\text{ is
bounded}\label{bu}%
\end{equation}
Moreover, by Lemma \ref{pis} , $u_{0}\geq0.$ Then the sequence $\left|
\nabla\theta\right|  ^{2}\ u_{0}\varphi_{n}$ is monotone and it converges
a.e$.$ to $\left|  \nabla\theta\right|  ^{2}\ u_{0}v^{+}.$ Then, by the
monotone convergence theorem, we get
\begin{equation}
\int\left|  k\nabla\theta\right|  ^{2}\ u_{0}\varphi_{n}dx\rightarrow
\int\left|  k\nabla\theta\right|  ^{2}\ u_{0}v^{+}dx\label{bubu}%
\end{equation}
By (\ref{bu}) and (\ref{bubu}) we get (\ref{fin}).

Taking the limit in (\ref{simbol}) and by using (\ref{pa}), (\ref{bi}),
(\ref{ba}), (\ref{fin}) we have
\[
\int\nabla u_{0}\cdot\nabla v^{+}+\left[  \left|  \mathbf{A}_{0}-k\nabla
\theta\right|  ^{2}-\left(  \phi_{0}-\omega\right)  ^{2}\right]  \,u_{0}%
v^{+}+W^{\prime}\left(  u_{0}\right)  v^{+}=0\
\]

Taking $\varphi_{n}=v^{-}\chi_{n}$ and arguing in the same way as before, we
get
\[
\int\nabla u_{0}\cdot\nabla v^{-}+\left[  \left\vert \mathbf{A}_{0}%
-k\nabla\theta\right\vert ^{2}-\left(  \phi_{0}-\omega\right)  ^{2}\right]
\,u_{0}v^{-}+W^{\prime}\left(  u_{0}\right)  v^{-}=0.
\]
Then
\[
\int\nabla u_{0}\cdot\nabla v+\left[  \left\vert \mathbf{A}_{0}-k\nabla
\theta\right\vert ^{2}-\left(  \phi_{0}-\omega\right)  ^{2}\right]
\,u_{0}v+W^{\prime}\left(  u_{0}\right)  v=0.
\]
Since $v\in\mathcal{D}$ is arbitrary, we get that equation (\ref{im}) is solved.
\end{proof}

\section{The natural constraints}

The functional $J$ presents two main difficulties:

\bigskip

1) The term $\int\left|  \nabla\times\mathbf{A}\right|  ^{2}$ is not a Sobolev
norm and it does not yield a control on $\int\left|  \nabla\mathbf{A}\right|
^{2}=\left\|  \mathbf{A}\right\|  _{\left(  \mathcal{D}^{1,2}\right)  ^{3}%
}^{2}.$

\bigskip

2) The presence of the term $-\int\left\vert \nabla\phi\right\vert ^{2}$ gives
to the functional $J$ a strong\ indefiniteness, namely any critical point of
$J$ has infinite Morse index: this fact is a great obstacle to a direct study
of its critical points.

\bigskip

In order to avoid these difficulties we introduce a suitable manifold
$\mathcal{M}\subset H$ such that:

\begin{itemize}
\item the critical points of $J$ restricted to $\mathcal{M}$ satisfy Eq.
(\ref{z1}), (\ref{z3}) ,\ref{z4}); namely $\mathcal{M}$ is a ''natural
constraint''\ for $J$.

\item The components $\mathbf{A}$ of the elements in $\mathcal{M}$ are
divergence free, then the term $\int\left|  \nabla\times\mathbf{A}\right|
^{2}$ can be replaced by $\left\|  \mathbf{A}\right\|  _{\left(
\mathcal{D}^{1,2}\right)  ^{3}}^{2}=\int\left|  \nabla\mathbf{A}\right|  ^{2}$.

\item The functional $I:=J\mid_{\mathcal{M}}$ does not exhibit the severe
indefiniteness before mentioned.
\end{itemize}

\subsection{The manifold of divergence free vector fields}

We shall denote by $u=u(r,x_{3})$ the maps having cylindrical symmetry, i.e.
those maps in $\mathbb{R}^{3}$ which depends only from $r=\sqrt{x_{1}%
^{2}+x_{2}^{2}}$ and $x_{3}.$ We set
\begin{equation}
\mathcal{A}_{0}:=\left\{  \mathbf{X}\in\mathcal{C}_{0}^{\infty}(\mathbb{R}%
^{3}\backslash\Sigma,\mathbb{R}^{3}):\mathbf{X}=b\left(  r,x_{3}\right)
\nabla\theta;\ b\in C_{0}^{\infty}\left(  \mathbb{R}^{3}\backslash
\Sigma,\mathbb{R}\right)  \right\} \label{ac}%
\end{equation}
Let $\mathcal{A}$ denote the closure of $\mathcal{A}_{0}$ with respect to the
norm of $\left(  \mathcal{D}^{1,2}\right)  ^{3}.$ We set
\begin{equation}
\mathcal{D}_{r}=\left\{  u\in\mathcal{D}:u=u(r,x_{3})\right\} \label{cf}%
\end{equation}
and we shall consider the following space
\begin{equation}
V:=\hat{H}_{r}^{1}\times\mathcal{A}\label{set}%
\end{equation}

where $\hat{H}_{r}^{1},$ is the closure of $\mathcal{D}_{r}$ with respect to
the $\hat{H}^{1}$ norm.

We need some technical preliminaries

\begin{lemma}
\label{contaiolo} If $\mathbf{A\in}\mathcal{A},\ $then $\nabla\times
\nabla\times\mathbf{A\in}\mathcal{A}^{\prime}$
\end{lemma}

\begin{proof}
Let $a\nabla\theta\in\mathcal{A}_{0},$ where $\ a\in C_{0}^{\infty}\left(
\mathbb{R}^{3}\backslash\Sigma,\mathbb{R}\right)  ,$ $a=a\left(
r,x_{3}\right)  .$ Then a straight computation shows that
\begin{equation}
\nabla\times\nabla\times(a\nabla\theta)=b\nabla\theta\label{king}%
\end{equation}

where
\[
b=-\frac{\partial^{2}a}{\partial r^{2}}+\frac{1}{r}\frac{\partial a}{\partial
r}-\frac{\partial^{2}a}{\partial x_{3}^{2}}%
\]
Now let $a_{n}\nabla\theta$ be a sequence in $\mathcal{A}_{0}$ converging in
$\left(  \mathcal{D}^{1,2}\right)  ^{3}$ to $\mathbf{A\in}\mathcal{A}$. By
continuity, $\nabla\times\nabla\times(a_{n}\nabla\theta)$ converges in
($\left(  \mathcal{D}^{1,2}\right)  ^{3})^{\prime}$ to $\nabla\times
\nabla\times\mathbf{A.}$ On the other hand by (\ref{king}) we have
\[
\nabla\times\nabla\times(a_{n}\nabla\theta)=b_{n}\nabla\theta\in
\mathcal{A}_{0}%
\]
Then, by definition, $\nabla\times\nabla\times\mathbf{A}\in\mathcal{A}%
^{\prime}$
\end{proof}

\begin{lemma}
\label{co}If $\mathbf{A}\in\mathcal{A}$ and $u\in\hat{H}_{r}^{1},$
\begin{equation}
\left(  \mathbf{A}-\nabla\theta\right)  u^{2}\in\mathcal{A}^{\prime
}\label{sat}%
\end{equation}

\end{lemma}

\begin{proof}
Set
\[
\mathcal{D}_{r}\left(  \mathbb{R}^{3}\backslash\Sigma\right)  =\left\{  u\in
C_{0}^{\infty}(\mathbb{R}^{3}\backslash\Sigma):u=u(r,x_{3})\right\}
\]
Since $u\in\hat{H}_{r}^{1}$, there exists a sequence $\left\{  u_{n}\right\}
\subset\mathcal{D}_{r}$ $\left(  \mathbb{R}^{3}\backslash\Sigma\right)  $ such
that
\begin{equation}
u_{n}\rightarrow u\text{ in }\hat{H}_{r}^{1}\label{conv}%
\end{equation}

Since $\mathbf{A\in}\mathcal{A}$, there exists a sequence $\left\{
b_{n}\right\}  \subset\mathcal{D}_{r}\left(  \mathbb{R}^{3}\backslash
\Sigma\right)  $ such that
\begin{equation}
b_{n}\nabla\theta\rightarrow\mathbf{A}\text{ in }(\mathcal{D}_{r}^{1,2}%
)^{3}\label{coco}%
\end{equation}

By (\ref{coco}) and (\ref{conv}) we deduce, following analogous arguments as
those used in proving lemma \ref{tecnico}, that
\begin{equation}
(b_{n}-1)u_{n}^{2}\nabla\theta=\left(  b_{n}\nabla\theta-\nabla\theta\right)
u_{n}^{2}\rightarrow\left(  \mathbf{A}-\nabla\theta\right)  u^{2}\text{ in
}\left(  (\mathcal{D}_{r}^{1,2})^{3}\right)  ^{\prime}\label{tt11}%
\end{equation}

Clearly
\begin{equation}
(b_{n}-1)u_{n}^{2}\in\mathcal{D}_{r}\left(  \mathbb{R}^{3}\backslash
\Sigma\right) \label{tt2}%
\end{equation}

Then, by (\ref{tt11}), (\ref{tt2}), we get (\ref{sat}).
\end{proof}

\bigskip

\begin{lemma}
\label{palla} If $\mathbf{A}\in\mathcal{A}$, then
\[
\int\left|  \nabla\times\mathbf{A}\right|  ^{2}=\int\left|  \nabla
\mathbf{A}\right|  ^{2}%
\]
and hence
\[
-\Delta\mathbf{A}=\nabla\times\nabla\times\mathbf{A}%
\]

\end{lemma}

\begin{proof}
Let $\mathbf{A=}b\nabla\theta\in\mathcal{A}_{0}.$ Since $b\ $depends only on
$r\ $and $x_{3},$\ it\ is easy to check that
\[
\nabla b\cdot\nabla\theta=0
\]
Since $\theta$ is harmonic in $\mathbb{R}^{3}\backslash\Sigma$ and $b$ has
support in $\mathbb{R}^{3}\backslash\Sigma$%
\[
b\Delta\theta=0
\]
Then
\[
\nabla\cdot\mathbf{A=}\nabla\cdot\left(  b\nabla\theta\right)  =\nabla
b\cdot\nabla\theta+b\Delta\theta=0
\]
Thus, by continuity, we get
\[
\int\left(  \nabla\cdot\mathbf{A}\right)  ^{2}=0\text{ for any }\mathbf{A}%
\in\mathcal{A}\text{ }%
\]
Then
\[
\int\left|  \nabla\times\mathbf{A}\right|  ^{2}=\int\left(  \nabla
\cdot\mathbf{A}\right)  ^{2}+\int\left|  \nabla\times\mathbf{A}\right|
^{2}=\int\left|  \nabla\mathbf{A}\right|  ^{2}%
\]
and clearly
\[
-\Delta\mathbf{A}=\nabla\times\nabla\times\mathbf{A}%
\]

\end{proof}

\bigskip

\subsection{The Gauss equation}

Equation (\ref{z3}), which we call Gauss equation, can be written as follows
\begin{equation}
-\Delta\phi+u^{2}\phi=\omega u^{2}\label{a2bis}%
\end{equation}

\begin{lemma}
\label{tre}Let $u\in H^{1}(\mathbb{R}^{3}),$ then there exists a unique
solution $\phi\in\mathcal{D}^{1,2}$ of (\ref{a2bis}).
\end{lemma}

\begin{proof}
$H^{1}(\mathbb{R}^{3})$ is continuously embedded into $L^{6}(\mathbb{R}^{3}),
$ then clearly
\begin{equation}
u^{2}\in L^{1}(\mathbb{R}^{3})\cap L^{3}(\mathbb{R}^{3})\label{interpol}%
\end{equation}
and, by interpolation, we have
\begin{equation}
u^{2}\in L^{\frac{3}{2}}(\mathbb{R}^{3})\label{Sob}%
\end{equation}

Now consider the bilinear form
\[
a(\phi,v)=\int\left\{  (\nabla\phi\mid\nabla v)+u^{2}\phi v\right\}  dx,\text{
}u,\phi\in\mathcal{D}^{1,2}%
\]
By (\ref{Sob}) we easily derive that $a(\phi,v)$ is equivalent to the standard
inner product in $\mathcal{D}^{1,2},$ i.e$.$
\[
c_{1}\left\|  \phi\right\|  _{\mathcal{D}^{1,2}}^{2}\leq a(\phi,\phi)\leq
c_{2}\left\|  \phi\right\|  _{\mathcal{D}^{1,2}}^{2},\text{ }c_{1},c_{2}>0
\]
On the other hand, using again (\ref{interpol}), we have
\[
u^{2}\in L^{\frac{6}{5}}(\mathbb{R}^{3})\subset\left(  \mathcal{D}%
^{1,2}\right)  ^{\prime}%
\]
So there exists a unique $\phi\in\mathcal{D}^{1,2}$ such that
\[
a(\phi,v)=-\omega\int u^{2}v\ \text{ for all }v\in\mathcal{D}^{1,2}.
\]
$\phi$ clearly solves (\ref{a2bis}).
\end{proof}

\ 

Clearly, if $u\in\hat{H}_{r}^{1}(\mathbb{R}^{3}),$ the solution $\phi$
$=\phi_{u}$ of (\ref{a2bis}) belongs to $\mathcal{D}_{r}^{1,2}.$ Then we can
define the map
\begin{equation}
u\in\hat{H}_{r}^{1}(\mathbb{R}^{3})\rightarrow\ Z\left(  u\right)  =\phi
_{u}\in\mathcal{D}_{r}^{1,2}\text{ solution of (\ref{a2bis})}\label{map}%
\end{equation}
Standard arguments show that the map $Z$ is $C^{1}.$ Since $\phi_{u}$ solves
(\ref{a2bis}), clearly we have
\begin{equation}
d_{\phi}J(u,\phi_{u},\mathbf{A)}=0\label{by}%
\end{equation}
For $u\in H^{1}(\mathbb{R}^{3}),$ let $\Phi=\Phi_{u}$ be the solution of the
equation (\ref{a2bis}) with $\omega=1,$ then $\Phi_{u}$ solves the equation
\begin{equation}
-\Delta\Phi_{u}+u^{2}\Phi_{u}=u^{2}\label{b2}%
\end{equation}
Clearly
\begin{equation}
\phi_{u}=\omega\Phi_{u}\label{d}%
\end{equation}
\ 

\begin{lemma}
\label{massimo}For any $u\in H^{1}(\mathbb{R}^{3}),$%
\[
0\leq\Phi_{u}\leq1
\]

\end{lemma}

\begin{proof}
By the \ maximum principle
\[
\Phi_{u}\geq0
\]
Moreover, arguing by contradiction assume that there is a region $\Omega$
where $\Phi_{u}>1$ and $\Phi_{u}=1$ on $\partial\Omega.$ Clearly
\[
-\Delta\left(  \Phi_{u}-1\right)  +u^{2}\left(  \Phi_{u}-1\right)
=-\Delta\Phi_{u}+u^{2}\Phi_{u}-u^{2}=0
\]
Then $v=\Phi_{u}-1$ solves the Dirichlet problem
\[
-\Delta v+u^{2}v=0\text{ in }\Omega,\text{ }v=0\text{ on }\partial\Omega
\]
Multiplying by $v$ and integrating in $\Omega$ we get
\[
\int_{\Omega}\left(  \left|  \nabla v\right|  ^{2}+u^{2}v^{2}\right)  dx=0
\]
Then $v=\Phi_{u}-1=0$ in $\Omega$ contradicting $\Phi_{u}>1$ in $\Omega$.
\end{proof}

\subsection{The reduced functional}

Now, if $\left(  u,\mathbf{A}\right)  \in\hat{H}_{r}^{1}\times\mathcal{A},$ we
set
\[
I(u,\mathbf{A})=J(u,Z\left(  u\right)  ,\mathbf{A})
\]
where $J$ is defined in (\ref{functional}). Thus $I\ $can be regarded as $J$
restricted to the manifold
\begin{equation}
\mathcal{M=}\left\{  (u,Z\left(  u\right)  ,\mathbf{A}):(u,\mathbf{A})\in
\hat{H}_{r}^{1}\times\mathcal{A}\right\} \label{mm}%
\end{equation}
We will refer to $I(u,\mathbf{A})$ as the \textit{reduced functional}.

From (\ref{b2}) we have
\begin{equation}
\int u^{2}\Phi_{u}dx=\int\left\vert \nabla\Phi_{u}\right\vert ^{2}dx+\int
u^{2}\Phi_{u}^{2}dx\label{base}%
\end{equation}

Now, using lemma \ref{palla} and Eq. (\ref{d}), (\ref{base}), we have:
\begin{align*}
I(u,\mathbf{A})  & =\frac{1}{2}\int\left\vert \nabla u\right\vert
^{2}-\left\vert \nabla\phi_{u}\right\vert ^{2}+\left\vert \nabla
\times\mathbf{A}\right\vert ^{2}\\
& +\frac{1}{2}\int\left[  \left\vert \mathbf{A}-k\nabla\theta\right\vert
^{2}-\left(  \phi_{u}-\omega\right)  ^{2}\right]  \,u^{2}+\int W\left(
u\right) \\
& =\frac{1}{2}\int\left(  \left\vert \nabla u\right\vert ^{2}+\left\vert
\nabla\mathbf{A}\right\vert ^{2}+\left\vert \mathbf{A}-k\nabla\theta
\right\vert ^{2}u^{2}\right) \\
& -\frac{1}{2}\omega^{2}\int\left(  \left\vert \nabla\Phi_{u}\right\vert
^{2}+u^{2}\Phi_{u}^{2}+\,u^{2}-2u^{2}\Phi_{u}\right)  +\int W\left(  u\right)
\\
& =\frac{1}{2}\int\left\vert \nabla u\right\vert ^{2}+\left\vert
\nabla\mathbf{A}\right\vert ^{2}+\left\vert \mathbf{A}-k\nabla\theta
\right\vert ^{2}u^{2}\\
& +\frac{1}{2}\int(1-\omega^{2}\left[  1-{\Phi}_{u}\right]  )\,u^{2}-\int
F\left(  u\right)
\end{align*}

\ Since ${\Phi}_{u}\leq1$ and $\omega^{2}<1$, the functional $I$ contains only
one negative term ( $-\int F\left(  u\right)  $)$.$ Moreover $\int\left\vert
\nabla\mathbf{A}\right\vert ^{2}$ replaces the term $\int\left\vert
\nabla\times\mathbf{A}\right\vert ^{2}.$ As a conseguence the functional $I$
does not exhibit the difficulties mentioned at the beginning of this section.

It can be shown also that $\mathcal{M}$ defined in (\ref{mm}) is a natural
constraint for the criticizing sequences of $J$, namely the following theorem holds

\begin{theorem}
\label{finale+}Let $\left(  u_{n},\ \mathbf{A}_{n}\right)  \in\hat{H}_{r}%
^{1}\times\mathcal{A}$ be a sequence such that $\forall(v,\mathbf{V})\in
\hat{H}_{r}^{1}\times\mathcal{A}$
\begin{equation}
dI\left(  u_{n},\ \mathbf{A}_{n}\right)  \left[  v,\mathbf{V}\right]
\rightarrow0\ \label{1}%
\end{equation}
Then, we have that
\begin{align}
d_{u}J\left(  u_{n},Z\left(  u_{n}\right)  ,\mathbf{A}_{n}\right)  \left[
v\right]   & \rightarrow0\text{ for all }v\in\hat{H}^{1}\label{11}\\
d_{\phi}J\left(  u_{n},Z\left(  u_{n}\right)  ,\mathbf{A}_{n}\right)  \left[
w\right]   & =0\text{ for all }n\in\mathbb{N}\text{ and }w\in D^{1,2}%
\label{12}\\
d_{\mathbf{A}}J\left(  u_{n},Z\left(  u_{n}\right)  ,\mathbf{A}_{n}\right)
\left[  \mathbf{V}\right]   & \rightarrow0\text{ for all }\mathbf{V\in}\left(
D^{1,2}\right)  ^{3}\label{13}%
\end{align}
where $d_{u},$ $d_{\phi},$ $d_{\mathbf{A}}$ denote the partial differentials
with respect to $u,$ $\phi,$ $\mathbf{A.}$
\end{theorem}

\begin{proof}
Using the chain rule, we have for all $(v,\mathbf{V})\in\hat{H}_{r}^{1}%
\times\mathcal{A}$%
\begin{align*}
dI\left(  u_{n},\ \mathbf{A}_{n}\right)  \left[  v,\mathbf{V}\right]   &
=d_{u}I\left(  u_{n},\mathbf{A}_{n}\right)  \left[  v\right]  +d_{\mathbf{A}%
}I\left(  u_{n},\mathbf{A}_{n}\right)  \left[  \mathbf{V}\right] \\
& =d_{u}J\left(  u_{n},Z\left(  u_{n}\right)  ,\mathbf{A}_{n}\right)  \left[
v\right] \\
& +d_{\phi}J\left(  u_{n},Z\left(  u_{n}\right)  ,\mathbf{A}_{n}\right)
\left[  d_{u}Z\left(  u_{n}\right)  \left[  w\right]  \right] \\
& +d_{\mathbf{A}}J\left(  u_{n},Z\left(  u_{n}\right)  ,\mathbf{A}_{n}\right)
\left[  \mathbf{V}\right]
\end{align*}

By equation (\ref{by}),
\begin{equation}
d_{\phi}J\left(  u_{n},Z\left(  u_{n}\right)  ,\mathbf{A}_{n}\right)
=0\label{3}%
\end{equation}
then for all $(v,\mathbf{V})\in\hat{H}_{r}^{1}\times\mathcal{A}$%
\[
dI\left(  u_{n},\ \mathbf{A}_{n}\right)  \left[  v,\mathbf{V}\right]
=d_{u}J\left(  u_{n},Z\left(  u_{n}\right)  ,\mathbf{A}_{n}\right)  \left[
v\right]  +d_{\mathbf{A}}J\left(  u_{n},Z\left(  u_{n}\right)  ,\mathbf{A}%
_{n}\right)  \left[  \mathbf{V}\right]
\]

So by (\ref{1}) we have for all $(v,\mathbf{V})\in\hat{H}_{r}^{1}%
\times\mathcal{A}$
\begin{equation}
d_{u}J\left(  u_{n},Z\left(  u_{n}\right)  ,\mathbf{A}_{n}\right)  \left[
v\right]  \rightarrow0\label{3b}%
\end{equation}

\begin{equation}
d_{\mathbf{A}}J\left(  u_{n},Z\left(  u_{n}\right)  ,\mathbf{A}_{n}\right)
\left[  \mathbf{V}\right]  \rightarrow0\label{3c}%
\end{equation}

Equations (\ref{3}), (\ref{3b}), (\ref{3c}), written explicitly take the
following form
\[
\int\nabla u_{n}\cdot\nabla v+\left[  \left|  \mathbf{A}_{n}-k\nabla
\theta\right|  ^{2}-\left(  \phi_{n}-\omega\right)  ^{2}\right]
\,u_{n}v+W^{\prime}\left(  u_{n}\right)  v=\left\langle \eta_{n}%
,v\right\rangle _{\hat{H}^{1}}%
\]
\[
\int\nabla\phi_{n}\cdot\nabla w+\left(  \phi_{n}-\omega\right)  u_{n}^{2}w=0
\]
\[
\int\nabla\mathbf{A}_{n}\cdot\nabla\mathbf{V}+\left(  \mathbf{A}_{n}%
-k\nabla\theta\right)  u_{n}^{2}\cdot\mathbf{V}=\left\langle \mathbf{\zeta
}_{n},\mathbf{V}\right\rangle _{\left(  D^{1,2}\right)  ^{3}}%
\]
where $(v,\mathbf{V})\in\hat{H}_{r}^{1}\times\mathcal{A}$, $w\in D^{1,2},$
$\phi_{n}=Z\left(  u_{n}\right)  ,$ $\eta_{n}\rightarrow0$ in $\left(  \hat
{H}^{1}\right)  ^{\prime}$ and $\mathbf{\zeta}_{n}\rightarrow0$ in $\left[
\left(  D^{1,2}\right)  ^{3}\right]  ^{\prime}.$

$\eta_{n}$ and $\mathbf{\zeta}_{n}$ have the following expression:
\begin{align}
\eta_{n}  & =-\Delta u_{n}+\left[  \left|  \mathbf{A}_{n}-k\nabla
\theta\right|  ^{2}-\left(  \phi_{n}-\omega\right)  ^{2}\right]
\,u_{n}+W^{\prime}\left(  u_{n}\right) \label{eb2}\\
\mathbf{\zeta}_{n}  & =-\Delta\mathbf{A}_{n}+\left(  \mathbf{A}_{n}%
-k\nabla\theta\right)  u_{n}^{2}\label{eb3}%
\end{align}

Since the Laplace operator preserves the cylindrical symmetry and $\left(
u_{n},\ \mathbf{A}_{n}\right)  \in\hat{H}_{r}^{1}\times\mathcal{A}$, we have
that $\eta_{n}\in\left(  \hat{H}_{r}^{1}\right)  ^{\prime}.$ Moreover, by
lemma (\ref{co}), we have that $\mathbf{\zeta}_{n}\in\mathcal{A}^{^{\prime}}$.

Now take any $(v,\mathbf{V})\in\hat{H}^{1}\times\left(  D^{1,2}\right)  ^{3}$
and we set
\[
v=v_{1}+v_{2}\text{ with }v_{1}\in\hat{H}_{r}^{1},.v_{2}\in\left(  \hat{H}%
_{r}^{1}\right)  ^{\bot}%
\]
and
\[
\mathbf{V=V}_{1}+\mathbf{V}_{2}\text{ with }\mathbf{V}_{1}\in\mathcal{A}%
,\mathbf{V}_{2}\in\mathcal{A}^{\bot}%
\]
Then, since $\eta_{n}\in\left(  \hat{H}_{r}^{1}\right)  ^{\prime}$ and
$v_{2}\in\left(  \hat{H}_{r}^{1}\right)  ^{\bot},$ we have
\[
d_{u}J\left(  u_{n},Z\left(  u_{n}\right)  ,\mathbf{A}_{n}\right)  \left[
v\right]  =\left\langle \eta_{n},v_{1}\right\rangle _{\hat{H}^{1}%
}+\left\langle \eta_{n},v_{2}\right\rangle _{\hat{H}^{1}}=\left\langle
\eta_{n},v_{1}\right\rangle _{\hat{H}^{1}}%
\]
So, since $\left\langle \eta_{n},v_{1}\right\rangle _{\hat{H}^{1}}%
\rightarrow0,$ we get
\[
d_{u}J\left(  u_{n},Z\left(  u_{n}\right)  ,\mathbf{A}_{n}\right)  \left[
v\right]  \rightarrow0
\]
Analogously it can be shown that
\[
d_{\mathbf{A}}J\left(  u_{n},Z\left(  u_{n}\right)  ,\mathbf{A}_{n}\right)
\left[  \mathbf{V}\right]  \rightarrow0
\]

\end{proof}

\section{The existence proof}

For the existence proof we need to construct suitable Palais-Smale sequences
for the functional
\begin{equation}
I(u,\mathbf{A})=\frac{1}{2}\int\left\vert \nabla u\right\vert ^{2}+\left\vert
\nabla\mathbf{A}\right\vert ^{2}+\left\vert \mathbf{A}-k\nabla\theta
\right\vert ^{2}\,u^{2}-\frac{\omega^{2}}{2}\int\left[  1-\Phi_{u}\right]
\,u^{2}+\int W\left(  u\right) \label{buono}%
\end{equation}
in the space $V:=\hat{H}_{r}^{1}\times\mathcal{A}.$

To simplify the notations, we set
\[
U=\left(  u,\mathbf{A}\right)
\]
\[
\left\Vert U\right\Vert _{V}=\left\Vert u\right\Vert _{\hat{H}^{1}}+\left\Vert
\mathbf{A}\right\Vert _{\left(  \mathcal{D}^{1,2}\right)  ^{3}}%
\]

First we prove the existence of a Palais-Smale sequence for $I$, namely the
following lemma holds

\begin{lemma}
\label{stand}Assume that the function $W$ \ satisfies assumptions (\ref{n})...
(\ref{nnnn}) and that $\omega^{2}<1$. Then there exits a PS sequence at some
level $c>0,$ namely a sequence $\left\{  U_{n}\right\}  \subset V$ such that
\[
I\left(  U_{n}\right)  \rightarrow c>0\text{ and }dI\left(  U_{n}\right)
\rightarrow0\text{ in }V^{\prime}%
\]

\end{lemma}

\begin{proof}
By a variant of the well known mountain pass theorem \cite{ar}, the existence
of a PS sequence will be guaranteed if we show that the functional
(\ref{buono}) has the mountain pass geometry, namely if there exist
$\alpha,\rho>0$ and $U_{0}\in V$ with $\left\|  U_{0}\right\|  _{V}>\rho$,
such that:
\begin{align}
I(U)  & \geq\alpha\text{ for }\left\|  U\right\|  _{V}=\rho\text{ }%
\label{mpg}\\
\text{and }I(U_{0})  & \leq0\label{mpgg}%
\end{align}

Let us first prove (\ref{mpg}). In the following $C_{1},..,C_{4}$ will denote
positive constants.

We have that
\begin{align}
\int\left|  \mathbf{A}-k\nabla\theta\right|  ^{2}\,u^{2}  & \geq\int\left(
\left|  \mathbf{A}\right|  ^{2}-2\frac{\left|  \mathbf{A}\right|  }{r}%
+\frac{1}{r^{2}}\right)  \,u^{2}\nonumber\\
& =\int\left(  \left|  \mathbf{A}\right|  ^{2}-2\left(  \left|  \mathbf{A}%
\right|  \sqrt{2}\frac{1}{\sqrt{2}r}\right)  +\frac{1}{r^{2}}\right)
\,u^{2}\nonumber\\
& \geq\int\left(  \left|  \mathbf{A}\right|  ^{2}-\left(  \left|
\mathbf{A}\right|  \sqrt{2}\right)  ^{2}-\left(  \frac{1}{\sqrt{2}r}\right)
^{2}+\frac{1}{r^{2}}\right)  \,u^{2}\nonumber\\
& =\int\left(  \left|  \mathbf{A}\right|  ^{2}-2\left|  \mathbf{A}\right|
^{2}-\frac{1}{2r^{2}}+\frac{1}{r^{2}}\right)  u^{2}\nonumber\\
& =\int-\left|  \mathbf{A}\right|  ^{2}u^{2}+\frac{u^{2}}{2r^{2}}\nonumber\\
& \geq\int\frac{u^{2}}{2r^{2}}-\left\|  \mathbf{A}\right\|  _{L^{6}}%
^{2}\left\|  u\right\|  _{L^{3}}^{2}\nonumber\\
& \geq\int\frac{u^{2}}{2r^{2}}-\frac{1}{2}\left\|  \mathbf{A}\right\|
_{L^{6}}^{4}-\frac{1}{2}\left\|  u\right\|  _{L^{3}}^{4}\nonumber\\
& \geq\frac{1}{2}\int\frac{u^{2}}{r^{2}}-\frac{C_{1}}{2}\left\|
\mathbf{A}\right\|  _{\mathcal{D}^{1,2}}^{4}-\frac{C_{2}}{2}\left\|
u\right\|  _{\hat{H}^{1}}^{4}\label{otto}%
\end{align}

Moreover by (\ref{nn}), (\ref{nnn}), we have that
\begin{equation}
\int F\left(  u\right)  \leq\frac{c}{p-1}\left\|  u\right\|  _{L^{p}}^{p}\leq
C_{3}\left\|  u\right\|  _{\hat{H}^{1}}^{p}\label{nove}%
\end{equation}

For any $U=\left(  u,\mathbf{A}\right)  \in V$ we have
\begin{align*}
I(U)  & =\frac{1}{2}\int\left|  \nabla u\right|  ^{2}+\left|  \nabla
\mathbf{A}\right|  ^{2}+\left|  \mathbf{A}-k\nabla\theta\right|  ^{2}%
\,u^{2}-\frac{\omega^{2}}{2}\int\left[  1-\Phi_{u}\right]  \,u^{2}+\int
W\left(  u\right) \\
(\text{by (\ref{otto}) and since }\omega & \geq{\phi}_{u}\geq0\text{)}%
\geq\frac{1}{2}\int\left|  \nabla u\right|  ^{2}+\left|  \nabla\mathbf{A}%
\right|  ^{2}+\\
& \frac{1}{2}\int\frac{u^{2}}{r^{2}}-\frac{C_{1}}{2}\left\|  \mathbf{A}%
\right\|  _{\mathcal{D}^{1,2}}^{4}-\frac{C_{2}}{2}\left\|  u\right\|
_{\hat{H}^{1}}^{4}-\frac{\omega^{2}}{2}\int u^{2}+\int\left(  \frac{1}{2}%
u^{2}-F(u)\right) \\
(\text{by (}\ref{nove}))  & \geq\frac{1}{2}\left(  \left\|  \mathbf{A}%
\right\|  _{\mathcal{D}^{1,2}}^{2}-C_{1}\left\|  \mathbf{A}\right\|
_{\mathcal{D}^{1,2}}^{4}\right)  +\\
& \frac{1}{2}\int\left(  \left|  \nabla u\right|  ^{2}+\frac{u^{2}}{r^{2}%
}\right)  -\frac{C_{2}}{2}\left\|  u\right\|  _{\hat{H}^{1}}^{4}%
+\frac{1-\omega^{2}}{2}\int u^{2}-C_{3}\left\|  u\right\|  _{\hat{H}^{1}}%
^{p}\\
(\text{ since }\omega^{2}  & <1)\geq\frac{1}{2}\left(  \left\|  \mathbf{A}%
\right\|  _{\mathcal{D}^{1,2}}^{2}-C_{1}\left\|  \mathbf{A}\right\|
_{\mathcal{D}^{1,2}}^{4}\right)  +\\
& \frac{1-\omega^{2}}{2}\left\|  u\right\|  _{\hat{H}^{1}}^{2}-\frac{C_{2}}%
{2}\left\|  u\right\|  _{\hat{H}^{1}}^{4}-C_{3}\left\|  u\right\|  _{\hat
{H}^{1}}^{p}\\
& =\frac{1}{2}\left(  1-C_{1}\left\|  \mathbf{A}\right\|  _{\mathcal{D}^{1,2}%
}^{2}\right)  \left\|  \mathbf{A}\right\|  _{\mathcal{D}^{1,2}}^{2}+\\
& \left(  \frac{1-\omega^{2}}{2}-\frac{C_{2}}{2}\left\|  u\right\|  _{\hat
{H}^{1}}^{2}-C_{3}\left\|  u\right\|  _{\hat{H}^{1}}^{p-2}\right)  \left\|
u\right\|  _{\hat{H}^{1}}^{2}%
\end{align*}
Then, since $\omega^{2}<1,$ there exist $\alpha,\rho>0$ such that
\begin{equation}
I(U)\geq\alpha\text{ for }\left\|  U\right\|  _{V}=\rho\label{dieci}%
\end{equation}

Let us now prove (\ref{mpgg}). Take $u_{0}\in C_{0}^{\infty}(\mathbb{R}%
^{3}\backslash\Sigma)$ $u_{0}\geq0,\ u_{0}\neq0$. By (\ref{nnnn})
\[
F(s)\geq F(1)s^{p},\ s>0;
\]
then
\begin{equation}
\int F\left(  tu_{0}\right)  \geq\int F(1)\left|  tu_{0}\right|  ^{p}\geq
C_{4}t^{p}\label{undici}%
\end{equation}

Now we have
\[
I(tu_{0},0)=\frac{t^{2}}{2}\int\left(  \left|  \nabla u_{0}\right|  ^{2}%
+\frac{u_{0}^{2}}{r^{2}}\,\right)  dx-\frac{1}{2}\omega^{2}t^{2}\int\left[
1-{\Phi}_{tu_{0}}\right]  \,u_{0}^{2}+\int W\left(  tu_{0}\right)  \leq
\]
\[
(\text{since }1\geq{\Phi}_{tu_{0}})\leq\frac{t^{2}}{2}\int\left(  \left|
\nabla u_{0}\right|  ^{2}+\frac{u_{0}^{2}}{r^{2}}\,\right)  dx+\frac{t^{2}}%
{2}\int u_{0}^{2}-\int F\left(  tu_{0}\right)  \leq\text{ }%
\]
\[
(\text{by (\ref{undici}))}\leq C_{5}t^{2}-C_{4}t^{p}.
\]

Then (\ref{mpgg}) is satisfied if we take $U_{0}=(tu_{0},0)$ with $t$
sufficiently large. \ 
\end{proof}

\bigskip

Our functional is invariant for translations in the $x_{3}$-direction, namely
for $U\in V$ and $L\in\mathbb{R}$ we have
\[
I(T_{L}U)=I(U)
\]
where
\[
T_{L}\left(  U\right)  \left(  x_{1},x_{2},x_{3}\right)  =U\left(  x_{1}%
,x_{2},x_{3}+L\right)
\]
As consequence of this invariance we have that a PS sequence for the
functional $I$ does not contain in general a (strongly) convergent
subsequence. So, in order to prove the existence of non trivial critical
points of $I$, we shall use the following strategy:

(i) we prove that any PS sequence $\left(  u_{n},\ \mathbf{A}_{n}\right)  $ of
the functional $I$ is bounded in $\hat{H}_{r}^{1}\times\left(  \mathcal{D}%
^{1,2}\right)  ^{3}$ $.$

(ii) we prove that there exists a PS sequence $\left(  u_{n},\ \mathbf{A}%
_{n}\right)  $\ whose weak limit $\left(  u_{0},\ \mathbf{A}_{0}\right)  $
gives rise to a non trivial critical point $\left(  u_{0},\ \phi_{u_{0}%
},\mathbf{A}_{0}\right)  $ of $J.$

\begin{lemma}
\label{14}Assume that the function $W$ \ satisfies assumptions (\ref{n})...
(\ref{nnnn}) and assume that $\omega^{2}<\min(1,\frac{p-2}{2})$. Let $\left\{
U_{n}\right\}  =\left\{  \left(  u_{n},\ \mathbf{A}_{n}\right)  \right\}
\subset\hat{H}_{r}^{1}\times\mathcal{A}$ be a PS sequence for the functional
$I$. Then $\left\{  U_{n}\right\}  $ is bounded in $\hat{H}_{r}^{1}%
\times\left(  \mathcal{D}^{1,2}\right)  ^{3}$.
\end{lemma}

\begin{proof}
Since $\left\{  \left(  u_{n},\ \mathbf{A}_{n}\right)  \right\}  \subset V$ is
a PS sequence for $I$, we have
\begin{equation}
dI\left(  U_{n}\right)  =\eta_{n}\rightarrow0\text{ in }V^{\prime}\label{ps}%
\end{equation}

and
\begin{equation}
I\left(  U_{n}\right)  =c_{n}\rightarrow c\label{pps}%
\end{equation}

From (\ref{ps}) we get
\begin{equation}
\int\left\vert \nabla u_{n}\right\vert ^{2}+\left[  \left\vert \mathbf{A}%
_{n}-k\nabla\theta\right\vert ^{2}\,-\left(  \phi_{n}-\omega\right)
^{2}\right]  u_{n}^{2}+W^{\prime}\left(  u_{n}\right)  u_{n}=\left\langle
\eta_{n},u_{n}\right\rangle \label{one}%
\end{equation}

where
\[
\phi_{n}=\omega\Phi_{u_{n}}\text{ }%
\]

The expression (\ref{pps}) can be written
\begin{equation}
\frac{1}{2}\int\left\vert \nabla u_{n}\right\vert ^{2}+\left\vert
\nabla\mathbf{A}_{n}\right\vert ^{2}+\left\vert \mathbf{A}_{n}-k\nabla
\theta\right\vert ^{2}\,u_{n}^{2}-\frac{1}{2}\int\omega\left[  \omega-{\phi
}_{n}\right]  \,u_{n}^{2}+\int W\left(  u_{n}\right)  =c_{n}\label{twoo}%
\end{equation}

Multiplying (\ref{one}) by $\frac{1}{p}$ and subtracting from (\ref{twoo}) we have%

\begin{align}
& \left(  \frac{1}{2}-\frac{1}{p}\right)  \int\left(  \left\vert \nabla
u_{n}\right\vert ^{2}+\left\vert \mathbf{A}_{n}-k\nabla\theta\right\vert
^{2}\,u_{n}^{2}\right)  dx+\left(  \frac{1}{2}-\frac{1}{p}\right)
\int\left\vert \nabla\mathbf{A}_{n}\right\vert ^{2}+\nonumber\\
& \int\left(  W\left(  u_{n}\right)  -\frac{1}{p}W^{\prime}\left(
u_{n}\right)  u_{n}+\frac{1}{p}\left(  \phi_{n}-\omega\right)  ^{2}u_{n}%
^{2}-\frac{1}{2}\omega\left[  \omega-{\phi}_{n}\right]  \,u_{n}^{2}\right)
dx\nonumber\\
& =c_{n}-\frac{1}{p}\left\langle \eta_{n},u_{n}\right\rangle \label{fase1}%
\end{align}

Now
\begin{align}
& W\left(  u_{n}\right)  -\frac{1}{p}W^{\prime}\left(  u_{n}\right)
u_{n}+\frac{1}{p}\left(  \phi_{n}-\omega\right)  ^{2}u_{n}^{2}-\frac{1}%
{2}\omega\left[  \omega-{\phi}_{n}\right]  \,u_{n}^{2}=(\text{by }%
(\ref{n}))\nonumber\\
& =\left(  \frac{1}{2}-\frac{1}{p}\right)  u_{n}^{2}+\frac{1}{p}F^{\prime
}\left(  u_{n}\right)  u_{n}-F\left(  u_{n}\right)  +\nonumber\\
& \omega^{2}\left(  \frac{1}{p}-\frac{1}{2}\right)  u_{n}^{2}+\left(  \frac
{1}{2}-\frac{2}{p}\right)  \phi_{n}\omega u_{n}^{2}\\
& \geq(\text{by }(\ref{nnnn})\geq\left(  \frac{1}{2}-\frac{1}{p}\right)
u_{n}^{2}+\omega^{2}\left(  \frac{1}{p}-\frac{1}{2}\right)  u_{n}^{2}+\left(
\frac{1}{2}-\frac{2}{p}\right)  \phi_{n}\omega u_{n}^{2}\nonumber\\
& =\left(  \frac{1}{2}-\frac{1}{p}\right)  \left(  1-\omega^{2}\right)
u_{n}^{2}+\left(  \frac{1}{2}-\frac{2}{p}\right)  \phi_{n}\omega u_{n}%
^{2}\label{fase2}%
\end{align}

Consider first the case $p\geq4$ and assume $\omega^{2}<1.$

Then, since $\phi_{n}\omega=\Phi_{u_{n}}\omega^{2}\geq0,$ from (\ref{fase2})
we have that
\begin{equation}
W\left(  u_{n}\right)  -\frac{1}{p}W^{\prime}\left(  u_{n}\right)  u_{n}%
+\frac{1}{p}\left(  \phi_{n}-\omega\right)  ^{2}u_{n}^{2}-\frac{1}{2}%
\omega\left[  \omega-{\phi}_{n}\right]  \,u_{n}^{2}\geq\left(  \frac{1}%
{2}-\frac{1}{p}\right)  \left(  1-\omega^{2}\right)  u_{n}^{2}\label{fase3}%
\end{equation}
Since $\omega^{2}<1$, from (\ref{fase3}) we deduce that
\begin{equation}
\int\left(  W\left(  u_{n}\right)  -\frac{1}{p}W^{\prime}\left(  u_{n}\right)
u_{n}+\frac{1}{p}\left(  \phi_{n}-\omega\right)  ^{2}u_{n}^{2}-\frac{1}%
{2}\omega\left[  \omega-{\phi}_{n}\right]  \,u_{n}^{2}\right)  dx\geq
C_{1}\left\|  u_{n}\right\|  _{L^{2}}^{2}\label{fase4}%
\end{equation}
where
\[
C_{1}=\left(  \frac{1}{2}-\frac{1}{p}\right)  \left(  1-\omega^{2}\right)  >0
\]
By (\ref{fase1}) and (\ref{fase4}) we get
\[
\left(  \frac{1}{2}-\frac{1}{p}\right)  \int\left(  \left|  \nabla
u_{n}\right|  ^{2}+\left|  \mathbf{A}_{n}-k\nabla\theta\right|  ^{2}u_{n}%
^{2}\right)  +\frac{1}{2}\int\left|  \nabla\mathbf{A}_{n}\right|  ^{2}%
+C_{1}\left\|  u_{n}\right\|  _{L^{2}}^{2}\leq
\]
\[
\leq c_{n}-\frac{1}{p}\left\langle \eta_{n},u_{n}\right\rangle
\]
Then
\[
C_{2}\left\|  u_{n}\right\|  _{H^{1}}^{2}+\frac{1}{2}\left\|  \mathbf{A}%
_{n}\right\|  _{\left(  \mathcal{D}^{1,2}\right)  ^{3}}^{2}\leq c_{n}+\frac
{1}{p}\left\|  \eta_{n}\right\|  _{H^{-1}}\left\|  u_{n}\right\|  _{H^{1}}%
\]

So we conclude that the sequences $\left\{  \left\|  u_{n}\right\|  _{H^{1}%
}^{2}\right\}  $ and $\left\{  \left\|  \mathbf{A}_{n}\right\|  _{\left(
\mathcal{D}^{1,2}\right)  ^{3}}^{2}\right\}  $ are bounded.

Now consider the case $p<4$ and assume $\omega^{2}<\frac{p-2}{2}$

By (\ref{fase2}), we easily have%

\begin{align*}
& W\left(  u_{n}\right)  -\frac{1}{p}W^{\prime}\left(  u_{n}\right)
u_{n}+\frac{1}{p}\left(  \phi_{n}-\omega\right)  ^{2}u^{2}-\frac{1}{2}%
\omega\left[  \omega-{\phi}_{n}\right]  \,u^{2}\geq(\text{ since }\Phi_{u_{n}%
}\leq1)\\
& \geq\left[  \left(  \frac{1}{2}-\frac{1}{p}\right)  \left(  1-\omega
^{2}\right)  +\left(  \frac{1}{2}-\frac{2}{p}\right)  \omega^{2}\right]
u_{n}^{2}\geq C_{3}u_{n}^{2}%
\end{align*}
Since $\omega^{2}<\frac{p-2}{2},$ we have
\[
C_{3}=\left(  \frac{1}{2}-\frac{1}{p}\right)  \left(  1-\omega^{2}\right)
+\left(  \frac{1}{2}-\frac{2}{p}\right)  \omega^{2}>0
\]
$.$

Then
\begin{equation}
\int\left(  W\left(  u_{n}\right)  -\frac{1}{p}W^{\prime}\left(  u_{n}\right)
u_{n}+\frac{1}{p}\left(  \phi_{n}-\omega\right)  ^{2}u_{n}^{2}-\frac{1}%
{2}\omega\left[  \omega-{\phi}_{n}\right]  \,u_{n}^{2}\right)  \geq
C_{3}\left\|  u_{n}\right\|  _{L^{2}}^{2}\label{fase5}%
\end{equation}

Arguing as in the first case, by (\ref{fase1}) and (\ref{fase5}), we again
conclude that the sequences $\left\{  \left\|  u_{n}\right\|  _{H^{1}}%
^{2}\right\}  $ and $\left\{  \left\|  \mathbf{A}_{n}\right\|  _{\left(
\mathcal{D}^{1,2}\right)  ^{3}}^{2}\right\}  $ are bounded.

It remains to prove that $u_{n}$ is bounded also in $\hat{H}^{1}.$

We have that
\begin{align*}
c+1  & \geq I(u_{n},\mathbf{A}_{n})\\
& =\frac{1}{2}\int\left|  \nabla u_{n}\right|  ^{2}+\left|  \nabla
\mathbf{A}_{n}\right|  ^{2}+\left|  \mathbf{A}_{n}-k\nabla\theta\right|
^{2}u_{n}^{2}-\frac{\omega^{2}}{2}\int\left[  1-\Phi_{u_{n}}\right]
\,u_{n}^{2}+\int W\left(  u_{n}\right) \\
& \geq\frac{1}{2}\int\left|  \nabla u_{n}\right|  ^{2}+\left|  \nabla
\mathbf{A}_{n}\right|  ^{2}+\left|  \mathbf{A}_{n}-k\nabla\theta\right|
^{2}\,u_{n}^{2}+u_{n}^{2}-C_{4}\\
& \geq\frac{1}{2}\int\left|  \nabla u_{n}\right|  ^{2}+\left|  \nabla
\mathbf{A}_{n}\right|  ^{2}+\left|  \mathbf{A}_{n}\right|  ^{2}\,u_{n}%
^{2}+\left|  k\nabla\theta\right|  ^{2}\,u_{n}^{2}-2\left|  \mathbf{A}%
_{n}\right|  \left|  k\nabla\theta\right|  \,u_{n}^{2}+u_{n}^{2}-C_{4}\\
& \geq\frac{1}{2}\left\|  u_{n}\right\|  _{\hat{H}^{1}}^{2}-\int\frac{k}%
{r}\left|  \mathbf{A}_{n}\right|  u^{2}-C_{4}%
\end{align*}

Also, we have that
\begin{align*}
\int\frac{k}{r}\left|  \mathbf{A}_{n}\right|  \,u_{n}^{2}  & \leq\int\frac
{1}{2}\left(  4k^{2}\left|  \mathbf{A}_{n}\right|  ^{2}+\frac{1}{4r^{2}%
}\right)  \,u_{n}^{2}=2k^{2}\int\left|  \mathbf{A}_{n}\right|  ^{2}%
\,u^{2}+\frac{1}{8}\int\frac{1}{r^{2}}\,u_{n}^{2}\\
& \leq\frac{1}{8}\left\|  u_{n}\right\|  _{\hat{H}^{1}}^{2}+C_{5}%
\end{align*}
Then
\begin{align*}
c+1  & \geq\frac{1}{2}\left\|  u_{n}\right\|  _{\hat{H}^{1}}^{2}-\frac{1}%
{8}\left\|  u_{n}\right\|  _{\hat{H}^{1}}^{2}-C_{5}-C_{4}\\
& =\frac{3}{8}\left\|  u_{n}\right\|  _{\hat{H}^{1}}^{2}-C_{5}-C_{4}%
\end{align*}

\end{proof}

\bigskip

\begin{lemma}
\label{2}Let the assumptions of Lemma \ref{stand} be satisfied. Then there is
a PS sequence $U_{n}=\left(  u_{n},\ \mathbf{A}_{n}\right)  \subset V$ for the
functional $I$ such that for $n$ large enough
\begin{equation}
\left\Vert u_{n}\right\Vert _{L^{p}}^{p}\geq C_{5}>0\text{ }\label{kappi}%
\end{equation}

\end{lemma}

\begin{proof}
By Lemma \ref{stand} there is a\ sequence $U_{n}=\left(  u_{n},\ \mathbf{A}%
_{n}\right)  \subset V$ satisfying assumptions (\ref{ps}) and (\ref{pps}) with
$c>0$. Clearly we have
\begin{equation}
\int\left(  \left|  \nabla u_{n}\right|  ^{2}+\left[  \left|  \mathbf{A}%
_{n}-k\nabla\theta\right|  ^{2}\,-\left(  \phi_{n}-\omega\right)  ^{2}\right]
u_{n}^{2}+W^{\prime}\left(  u_{n}\right)  u_{n}\right)  dx=\left\langle
\eta_{n},u_{n}\right\rangle \label{compl}%
\end{equation}

where $\left\|  \eta_{n}\right\|  _{\hat{H}^{-1}}\rightarrow0.$

By Lemma \ref{14} the sequence $\left\{  \left\Vert u_{n}\right\Vert _{\hat
{H}^{1}}\right\}  $ is bounded. Then
\begin{equation}
\int\left(  \left\vert \nabla u_{n}\right\vert ^{2}+\left[  \left\vert
\mathbf{A}_{n}-k\nabla\theta\right\vert ^{2}\,-\left(  \phi_{n}-\omega\right)
^{2}\right]  u_{n}^{2}+W^{\prime}\left(  u_{n}\right)  u_{n}\right)
dx\rightarrow0\label{primula}%
\end{equation}

Moreover, by (\ref{n}), (\ref{nnn}) and (\ref{nnnn}), we get from
(\ref{compl})
\[
\left\|  u_{n}\right\|  _{H^{1}}^{2}+\int\left|  \mathbf{A}_{n}-k\nabla
\theta\right|  ^{2}u_{n}^{2}-\omega^{2}\int\left(  \Phi_{u_{n}}-1\right)
^{2}u_{n}^{2}\leq\varepsilon_{n}+\left\|  u_{n}\right\|  _{L^{p}}^{p},\text{
}\varepsilon_{n}\rightarrow0
\]
Then, since $0\leq\Phi_{u_{n}}\leq1,$ we have
\begin{equation}
\left(  1-\omega^{2}\right)  \left\|  u_{n}\right\|  _{H^{1}}^{2}%
\leq\varepsilon_{n}+\left\|  u_{n}\right\|  _{L^{p}}^{p}\label{kappa}%
\end{equation}
Clearly, since $\omega^{2}<1,$ (\ref{kappi}) will be a consequence of
(\ref{kappa}) if we show that
\begin{equation}
0<C_{3}\leq\left\|  u_{n}\right\|  _{H^{1}}^{2},\text{ }C_{3}>0\label{kappo}%
\end{equation}

We argue by contradiction and assume that (up to a subsequence)
\begin{equation}
\left\|  u_{n}\right\|  _{H^{1}}^{2}\rightarrow0\label{contr}%
\end{equation}

Since $\left\{  \left(  u_{n},\ \mathbf{A}_{n}\right)  \right\}  $ is a PS
sequence for the functional $I$, we have
\begin{align}
-\Delta\phi_{n}-\left(  \omega-\phi_{n}\right)  ~u_{n}^{2}  & =0,\text{ }%
\phi_{n}=Z(u_{n})\\
-\Delta\mathbf{A}_{n}-\left(  k\nabla\theta-\mathbf{A}_{n}\right)  u_{n}^{2}
& =\eta_{n}\rightarrow0\text{ in }\ \left(  \left(  \mathcal{D}^{1,2}\right)
^{3}\right)  ^{\prime}%
\end{align}
From these two equations we get
\begin{align}
\left\Vert \phi_{n}\right\Vert _{D^{1,2}}^{2}  & =\int\left(  \omega-\phi
_{n}\right)  ~u_{n}^{2}\phi_{n}\label{get1}\\
\left\Vert \mathbf{A}_{n}\right\Vert _{\left(  D^{1,2}\right)  ^{3}}^{2}  &
=\int\left(  \mathbf{A}_{n}-k\nabla\theta\right)  \cdot\mathbf{A}_{n}u_{n}%
^{2}+\delta_{n}\label{get2}%
\end{align}

where $\delta_{n}=\left\langle \eta_{n},\mathbf{A}_{n}\right\rangle .$ Since,
by Lemma \ref{14}, $\left\{  \left\|  \mathbf{A}_{n}\right\|  _{\left(
\mathcal{D}^{1,2}\right)  ^{3}}^{2}\right\}  $ is bounded, we have
\begin{equation}
\delta_{n}\rightarrow0\label{inf}%
\end{equation}

From (\ref{get1}) we have
\begin{equation}
\left\Vert \phi_{n}\right\Vert _{D^{1,2}}^{2}\leq\omega\int~u_{n}^{2}\phi
_{n}\leq\omega\left\Vert u_{n}\right\Vert _{L^{\frac{12}{5}}}^{2}\left\Vert
\phi_{n}\right\Vert _{L^{6}}\leq C\left\Vert u_{n}\right\Vert _{L^{\frac
{12}{5}}}^{2}\left\Vert \phi_{n}\right\Vert _{D^{1,2}}\label{contra}%
\end{equation}

By (\ref{contr}) and (\ref{contra}) we deduce
\begin{equation}
\left\|  \phi_{n}\right\|  _{D^{1,2}}\rightarrow0\text{ }\label{dedu}%
\end{equation}
By (\ref{dedu}) and (\ref{contr}) we easily have
\begin{equation}
\int\left(  \left|  \nabla u_{n}\right|  ^{2}-\left(  \phi_{n}-\omega\right)
^{2}u_{n}^{2}+W^{\prime}\left(  u_{n}\right)  u_{n}\right)  \rightarrow
0\label{eas}%
\end{equation}
then by (\ref{primula}) and (\ref{eas}) we get
\begin{equation}
\int\left|  \mathbf{A}_{n}-k\nabla\theta\right|  ^{2}\,u_{n}^{2}%
\rightarrow0\label{finna}%
\end{equation}

Let us show that we have also
\begin{equation}
\left\|  \mathbf{A}_{n}\right\|  _{\left(  D^{1,2}\right)  ^{3}}%
\rightarrow0\label{deda}%
\end{equation}

Now
\begin{align}
\left\vert \int\left(  \mathbf{A}_{n}-k\nabla\theta\right)  \cdot
\mathbf{A}_{n}u_{n}^{2}\right\vert  & \leq\int\left\vert \mathbf{A}%
_{n}-k\nabla\theta\right\vert \left\vert \mathbf{A}_{n}\right\vert u_{n}^{2}\\
& \leq\left(  \int\left\vert \mathbf{A}_{n}-k\nabla\theta\right\vert ^{2}%
u_{n}^{2}\right)  ^{\frac{1}{2}}\left(  \int\left\vert \mathbf{A}%
_{n}\right\vert ^{2}u_{n}^{2}\right)  ^{\frac{1}{2}}\nonumber\label{sos}%
\end{align}
Since the sequence $\left\{  \left\Vert \mathbf{A}_{n}\right\Vert _{\left(
\mathcal{D}^{1,2}\right)  ^{3}}\right\}  $ is bounded, we easily get from
(\ref{contr}) that
\begin{equation}
\int\left\vert \mathbf{A}_{n}\right\vert ^{2}u_{n}^{2}\rightarrow0\label{neu}%
\end{equation}
By (\ref{finna}) and (\ref{neu}), we deduce that
\begin{equation}
\int\left(  \mathbf{A}_{n}-k\nabla\theta\right)  \cdot\mathbf{A}_{n}u_{n}%
^{2}\rightarrow0\label{nev}%
\end{equation}
By (\ref{nev}) and (\ref{get2}) we get (\ref{deda}).

By (\ref{contr}) and (\ref{dedu}) we have
\begin{equation}
\int(1-\omega^{2}\left[  1-\Phi_{u_{n}}\right]  )\,u_{n}^{2}-F(u_{n}%
)\rightarrow0,\label{zero}%
\end{equation}
Then, by (\ref{zero}), (\ref{contr}), (\ref{deda}) and (\ref{finna}), we
deduce that
\begin{align*}
I\left(  u_{n},\ \mathbf{A}_{n}\right)   & =\frac{1}{2}\int\left|  \nabla
u_{n}\right|  ^{2}+\left|  \nabla\mathbf{A}_{n}\right|  ^{2}+\left|
\mathbf{A}_{n}-k\nabla\theta\right|  ^{2}u_{n}^{2}\\
+\frac{1}{2}\int(1-\omega^{2}\left[  1-\Phi_{u_{n}}\right]  )\,u_{n}^{2}-\int
F(u_{n})  & \rightarrow0
\end{align*}
This contradicts the assumption $I\left(  U_{n}\right)  \rightarrow c>0$. Then
(\ref{kappo}) holds.
\end{proof}

\bigskip

\begin{lemma}
\label{prec}Let the assumptions of Lemma \ref{stand} be satisfied. Then there
is a $PS$ sequence $U_{n}=\left(  u_{n},\mathbf{A}_{n}\right)  $ for the
functional $I$ such that $U_{n}\rightarrow U_{0}=\left(  u_{0},\mathbf{A}%
_{0}\right)  $ weakly in $\hat{H}^{1}\times\left(  \mathcal{D}^{1,2}\right)
^{3}$ and $u_{0}\neq0$
\end{lemma}

\begin{proof}
By Lemma \ref{2} we can take a PS sequence $U_{n}=\left(  u_{n},\ \mathbf{A}%
_{n}\right)  \subset V$ for the functional $I$ such that
\begin{equation}
\left\Vert u_{n}\right\Vert _{L^{p}}^{p}\geq C_{5}>0\text{ for }n\text{ large}%
\end{equation}
By Lemma \ref{14} the sequence $\left\{  U_{n}\right\}  $ is bounded in
$\hat{H}^{1}\times\left(  \mathcal{D}^{1,2}\right)  ^{3}.$ So we can assume
that $U_{n}\rightarrow U_{0}$ weakly in $\hat{H}^{1}\times\left(
\mathcal{D}^{1,2}\right)  ^{3}.$ If $U_{0}\neq0$ we have finished. If not we
will show that there is a sequence of integers $j_{n}$ such that (up to a
subsequence) $V_{n}:=T_{j_{n}}U_{n}\rightarrow U_{0}\neq0$ weakly in
$H^{1}\times\left(  \mathcal{D}^{1,2}\right)  ^{3}$.

We set
\[
\Omega_{j}=\left\{  \left(  x_{1},x_{2},x_{3}\right)  :j\leq x_{3}%
<j+1\right\}  \text{, }j\text{ integer}%
\]
We have for $n$ large
\begin{align*}
0  & <C_{5}\leq\left\|  u_{n}\right\|  _{L^{p}}^{p}=\sum_{j}\int_{\Omega_{j}%
}\left|  u_{n}\right|  ^{p}=\sum_{j}\left(  \int_{\Omega_{j}}\left|
u_{n}\right|  ^{p}\right)  ^{\left(  p-2\right)  /p}\cdot\left(  \int
_{\Omega_{j}}\left|  u_{n}\right|  ^{p}\right)  ^{2/p}\\
& \leq\sup_{j}\left\|  u_{n}\right\|  _{L^{p}\left(  \Omega_{j}\right)
}^{p-2}\sum_{j}\left(  \int_{\Omega_{j}}\left|  u_{n}\right|  ^{p}\right)
^{2/p}\leq C_{6}\cdot\sup_{j}\left\|  u_{n}\right\|  _{L^{p}\left(  \Omega
_{j}\right)  }^{p-2}\cdot\sum_{j}\left\|  u_{n}\right\|  _{H^{1}\left(
\Omega_{j}\right)  }^{2}\\
& \leq C_{6}\cdot\sup_{j}\left\|  u_{n}\right\|  _{L^{p}\left(  \Omega
_{j}\right)  }^{p-2}\cdot\left\|  u_{n}\right\|  _{H^{1}\left(  \mathbb{R}%
^{3}\right)  }^{2}\leq(\text{since }\left\|  u_{n}\right\|  _{H^{1}\left(
\mathbb{R}^{3}\right)  }^{2}\leq C_{7})\text{ }\\
& \leq C_{6}C_{7}\sup_{j}\left\|  u_{n}\right\|  _{L^{p}\left(  \Omega
_{j}\right)  }^{p-2}%
\end{align*}
Then, for $n$ large$,$ there exists an integer $j_{n}$ such that
\begin{equation}
\left\|  u_{n}\right\|  _{L^{p}\left(  \Omega_{j_{n}}\right)  }^{p-2}\geq
\frac{C_{4}}{2C_{5}C_{6}}:=C_{7}>0\label{paracula}%
\end{equation}
Now set
\[
\left(  u_{n}^{\prime},\mathbf{A}_{n}^{^{\prime}}\right)  =\ U_{n}^{\prime
}(x_{1},x_{2},x_{3})=U_{n}(x_{1},x_{2},x_{3}+j_{n})=T_{j_{n}}\left(
U_{n}\right)
\]
By lemma \ref{14} the sequence $u_{n}^{\prime}$ is bounded $\hat{H}^{1}\left(
\mathbb{R}^{3}\right)  ,$ then (up to a subsequence) it converges weakly to
$u_{0}\in\hat{H}^{1}\left(  \mathbb{R}^{3}\right)  .$ We want to show that
$u_{0}\neq0.$ Now, let $\varphi=\varphi\left(  x_{3}\right)  $ be a
nonnegative, $C^{\infty}$-function whose value is $1$ for $0<x_{3}<1$ and $0$
for $\left|  x_{3}\right|  >2.$ Then, the sequence $\varphi u_{n}^{\prime} $
is bounded in $H_{0}^{1}(\mathbb{R}^{2}\times(-2,2)),$ moreover $\varphi
u_{n}^{\prime}$ has cylindrical symmetry. Then, using the compactness result
proved in \cite{el}, we have
\[
\varphi u_{n}^{\prime}\rightarrow\chi\text{ strongly in }L^{p}(\mathbb{R}%
^{2}\times(-2,2))
\]
On the other hand
\begin{equation}
\varphi u_{n}^{\prime}\rightarrow\varphi u_{0}\text{ }a.e\text{.}\label{punto}%
\end{equation}

Then
\begin{equation}
\varphi u_{n}^{\prime}\rightarrow\varphi u_{0}\text{ strongly in }%
L^{p}(\mathbb{R}^{2}\times(-2,2))\label{virgola}%
\end{equation}

Moreover by (\ref{paracula})
\begin{equation}
\left\|  \varphi u_{n}^{\prime}\right\|  _{L^{p}\left(  \mathbb{R}^{2}%
\times(-2,2)\right)  }\geq\left\|  u_{n}^{\prime}\right\|  _{L^{p}\left(
\Omega_{0}\right)  }=\left\|  u_{n}\right\|  _{L^{p}\left(  \Omega_{j_{n}%
}\right)  }\geq C_{7}^{1/\left(  p-2\right)  }\label{pp}%
\end{equation}

Then$\ $by (\ref{virgola}) and (\ref{pp})
\[
\left\|  \varphi u_{0}\right\|  _{L^{p}\left(  \mathbb{R}^{2}\times
(-2,2)\right)  }\geq C_{7}^{1/\left(  p-2\right)  }>0.
\]
Thus we have that $u_{0}\neq0.$
\end{proof}

\bigskip Now we complete the proof Theorem \ref{main}.

By Lemma \ref{prec} there is a PS sequence $U_{n}=\left(  u_{n},\mathbf{A}%
_{n}\right)  $ for the functional $I$ such that $U_{n}\rightarrow
U_{0}=\left(  u_{0},\mathbf{A}_{0}\right)  $ weakly in $\hat{H}^{1}%
\times\left(  \mathcal{D}^{1,2}\right)  ^{3}$ and $u_{0}\neq0$. We shall show
that $\left(  u_{0},Z(u_{0}),\mathbf{A}_{0}\right)  $ is a critical point of
$J.$

Since $\left(  u_{n},\mathbf{A}_{n}\right)  $ is a PS sequence for $I$ \ we
have for any $\left(  v,\mathbf{V}\right)  \in\hat{H}_{r}^{1}\times
\mathcal{A}$
\begin{equation}
d_{u}I\left(  u_{n},\mathbf{A}_{n}\right)  \left[  v\right]  \rightarrow
0\label{wea}%
\end{equation}

\begin{equation}
d_{\mathbf{A}}I\left(  u_{n},\mathbf{A}_{n}\right)  \left[  \mathbf{V}\right]
\rightarrow0\label{weaa}%
\end{equation}
Then, since the manifold $\mathcal{M}$ defined in (\ref{mm}) is a natural
constraint for the Palais-Smale sequences (see Theorem \ref{finale+}), for any
$v,w\in C_{0}^{\infty}(\mathbb{R}^{3}\backslash\Sigma)$ and $\mathbf{V}%
\in\left(  C_{0}^{\infty}(\mathbb{R}^{3}\backslash\Sigma)\right)  ^{3}$, we
get
\begin{equation}
d_{u}J(u_{n},Z(u_{n}),\mathbf{A}_{n})\left[  v\right]  \rightarrow0\label{ri}%
\end{equation}
\begin{equation}
d_{\phi}J(u_{n},Z(u_{n}),\mathbf{A}_{n})\left[  w\right]  =0\label{ru}%
\end{equation}

\begin{equation}
d_{\mathbf{A}}J(u_{n},Z(u_{n}),\mathbf{A}_{n})\left[  \mathbf{V}\right]
\rightarrow0\label{ris}%
\end{equation}

Since $v,w$ and $\mathbf{V}$ have compact support in $\mathbb{R}^{3}$%
$\backslash$%
$\Sigma$ and $\left(  u_{n},\mathbf{A}_{n}\right)  $ $\rightarrow\left(
u_{0},\mathbf{A}_{0}\right)  $ weakly in $\hat{H}^{1}\times\left(
\mathcal{D}^{1,2}\right)  ^{3}$, by standard argument we can take the limit in
(\ref{ri}), (\ref{ru}), (\ref{ris}) and we get that $u_{0},Z(u_{0}%
),\mathbf{A}_{0}$ solve the equations
\begin{equation}
d_{u}J(u_{0},Z(u_{0}),\mathbf{A}_{0})\left[  v\right]  =0\label{pip}%
\end{equation}
\begin{equation}
d_{\phi}J(u_{0},Z(u_{0}),\mathbf{A}_{0})\left[  w\right]  =0\label{pap}%
\end{equation}

\begin{equation}
d_{\mathbf{A}}J(u_{0},Z(u_{0}),\mathbf{A}_{0})\left[  \mathbf{V}\right]
=0\label{pup}%
\end{equation}

Now by a density and continuity argument the test functions $v,w$ and
$\mathbf{V}$ in (\ref{pip}), (\ref{pap}) and (\ref{pup}) can be taken in
$\hat{H}^{1},$ $\mathcal{D}^{1,2}$ and $\left(  \mathcal{D}^{1,2}\right)
^{3}$ respectively. Then $\left(  u_{0},Z(u_{0}),\mathbf{A}_{0}\right)  $ is a
critical point of $J.$ Therefore, by using Theorem \ref{finale}, we get that
$u_{0},Z(u_{0}),\mathbf{A}_{0}$ solve (\ref{z1}), (\ref{z3}), (\ref{z4}) in
the sense of distributions in $\mathbb{R}^{3}.$ Finally, since $u_{0}\neq0,$
the conclusion follows (see Remark \ref{iv}).

\end{document}